\documentclass[11pt]{amsart}
\usepackage{texdraw}
\usepackage{amsmath, amssymb}
\usepackage[all]{xy}

\newcommand{\vs}{\vskip 12 pt}
\newcommand{\ra}{\rightarrow}

\newcommand{\PP}{\mathbb P}
\newcommand{\BC}{\mathbb C}
\newcommand{\BZ}{\mathbb Z}

\newcommand{\CB}{\mathcal B}
\newcommand{\CR}{\mathcal R}

\newcommand{\CF}{\mathcal F}
\newcommand{\CS}{\mathcal S}
\newcommand{\CD}{\mathcal D}
\newcommand{\CZ}{\mathcal Z}

\newcommand{\BS}{\mathbb S}

\newcommand{\ba}{\mathbf a}
\newcommand{\bb}{\mathbf b}

\newcommand{\bI}{\mathbf I}
\newcommand{\bi}{\mathbf i}

\newcommand{\fg}{\mathfrak{g}}
\newcommand{\fm}{\mathfrak{m}}
\newcommand{\fn}{\mathfrak{n}}
\newcommand{\fp}{\mathfrak{p}}
\newcommand{\fq}{\mathfrak{q}}
\newcommand{\fb}{\mathfrak{b}}

\newcommand{\fl}{\mathfrak{l}}
\newcommand{\fh}{\mathfrak{h}}

\newtheorem{prop}{Proposition}[section]
\newtheorem{thm}[prop]{Theorem}
\newtheorem{lem}[prop]{Lemma}

\theoremstyle{remark}
  \newtheorem{rk}[prop]{Remark}
\theoremstyle{definition}
 \newtheorem{exam}[prop]{Example}
 \newtheorem{defn}[prop]{Definition}

\newenvironment{pf}[1]{\noindent {\it{Proof.}} {#1}}{$ \square $ \vs}

\begin{document}

\title{Rigidity of smooth Schubert varieties in Hermitian symmetric spaces}

\author{Jaehyun Hong}

\begin{abstract}

In this paper we study the space $\mathcal Z_k(G/P, r[X_w])$ of
effective $k$-cycles $X$ in $G/P$ with the homology class equal to
an integral multiple of the homology class of Schubert variety
$X_w$ of type $w$. When $X_w$ is a proper linear subspace $\mathbb
P^k$ $(k<n)$ of a linear space $\mathbb P^n$ in $G/P \subset
\mathbb P(V)$, we know that $\mathcal Z_k(\mathbb P^n, r[\mathbb
P^k])$ is already complicated. We will show that for a smooth
Schubert variety $X_w$ in a Hermitian symmetric space, any
irreducible subvariety $X$ with the homology class $[X]=r[X_w]$,
$r\in \mathbb Z$ is again a Schubert variety of type $w$, unless
$X_w$ is a non-maximal linear space. In particular, any local
deformation of such a smooth Schubert variety in Hermitian
symmetric space $G/P$ is obtained by the action of the Lie group
$G$.

\end{abstract}


\maketitle


\section{Introduction}

In this paper we study the space $\mathcal Z_k(G/P, r[X_w])$ of
effective $k$-cycles $X$ in $G/P$ with the homology class equal to
an integral multiple of the homology class of a Schubert variety
$X_w$ of type $w$. For example, $\CZ_1(\PP^n, [\PP^1])$ consists
of all $\PP^1$'s in $\PP^n$ so every element in $\CZ_1(\PP^n,
[\PP^1])$ is obtained by acting $g \in SL(n+1)$ on a fixed
$\PP^1$. However, $\CZ_1(\PP^n, r[\PP^1]), r>1$ contains not only
the sum $\PP^1 + \cdots + \PP^1$ of $r$ $\PP^1$'s but also the sum
$C_1 + \cdots + C_s$ of curves $C_i$ of degree $r_i$ in $\PP^n$,
where $r_1 + \cdots + r_s=r$.

In general the cycle space $\CZ_k(M, z)$ for $z \in H_{2k}(M,
\mathbb Z)$ is complicated. But, as in the case of $\PP^1$ in
$\PP^n$, for certain subvarieties $X$ of a homogeneous space
$G/P$, the cycle space $Z_k(G/P, [X])$ is rather simple: $G$ acts
on $\CZ_k(G/P, [X])$ transitively. In particular, any local
deformation of $X$ in $G/P$ is given by the action of $G$. Even
more,  for certain Schubert varieties $X_w$ in a Hermitian
symmetric space $G/P$, the cycle space $\CZ_k(G/P, r[X_w])$
consists of the sums of $r$ Schubert varieties of type $w$ and is
just the $r$-symmetric product of the $G$-homogenous space
$\CZ_k(G/P, [X_w])$.

\vs{\bf Question.} For which Schubert variety $X_w$ in Hermitian
symmetric space $G/P$, does the cycle space $\CZ_k(G/P, r[X_w])$
consist of the sums of $r$ Schubert varieties of type $w$? \vs

Walters and Bryant studied this problem by transforming it to the
problem on the integral varieties of differential
systems(\cite{W}, \cite{B}). From the fact that there exists a
closed positive differential form $\phi_w$ for  $w\in W^P$ such
that $\int_{X_{w}}\phi_v=0$ for  all $v \not= w$ with
$\ell(v)=\ell(w)$(\cite{K2}), it follows that

\begin{eqnarray*}
[X]=r[X_w]
 &\Leftrightarrow& \int_X
\phi_{v}=0, \forall v \not=w, \ell(v)=\ell(w) \\
 &\Leftrightarrow& \phi_{v}|_{X}=0, \forall v \not=w, \ell(v)=\ell(w).
 \end{eqnarray*}

\vs
 \noindent So such an $X$ satisfies a first order holomorphic  partial differential
 equation. In particular, if  $X$ is tangent to a Schubert variety of type $w$ at each point
$x \in X$, which may depend on the point $x$, then we have $[X]=
r[X_w]$ for some $r \in \mathbb Z$. Since the ray generated by the
Schubert cycle $X_w$ is extremal, if an effective cycle $X=X_1 +
\cdots + X_n$ has the homology class $[X]=r[X_w]$, then $r=r_1+
\cdots + r_n$ and $[X_i]=r_i[X_w]$, where $X_i$ is an irreducible
compact complex variety of dimension $k$. So we may assume that
$X$ is irreducible.

Define a differential system $\CR_w$ by putting together all the
tangent $\ell(w)$-subspaces $W$ such that $\phi_v|_W=0$ for all $v
\not=w$ with $\ell(v)=\ell(w)$. We say that the Schubert variety
$X_w$ is {\it Schur rigid} if any irreducible integral variety of
$\CR_{w}$ is a Schubert variety of type $w$, or equivalently,
$\CZ_k(M, r[X_w])$ consists only of the sums of  Schubert cycles
of type $w$. Now the question becomes how to compute $\CR_w$ and
how to find all the integral varieties of $\CR_w$.

The differential system $\CR_w$ contains the differential system
$\CB_w$ with the fiber consisting of the tangent space of all the
Schubert varieties $X_w$ of type $w$ passing through the given
point. We say that the Schubert variety $X_w$ is {\it Schubert
rigid} if any irreducible integral variety of $\CB_w$ is a
Schubert variety of type $w$. Thus  $X_w$ is Schur rigid  if
$\CB_w$ is equal to $\CR_{w}$ and $X_w$ is Schubert
rigid(\cite{W}, Section 2.8 of \cite{B}). Conversely, $X_w$ is not
Schur rigid if either $\CB_w$ is a proper subvariety of $\CR_w$ or
$\CB_w$ is not Schubert rigid(Proposition 2 of \cite{B},
Proposition \ref{equivalence}).

It is known that smooth Schubert varieties in Grassmannian
$Gr(m,n)$ are Schubert rigid unless it is a non-maximal linear
spaces in $Gr(m,n)$ (\cite{W}). Schur and Schubert rigidity of
several kinds of smooth and singular Schubert varieties  in
Hermitian symmetric spaces are investigated systematically  in
\cite{B}. The Schubert rigidity of linear spaces in general
homogeneous spaces with $b_2=1$ is studied in \cite{CH}.

\vs

In this paper we restrict ourselves to the case of  smooth
Schubert varieties $X_w$ in Hermitian symmetric spaces $G/P$.
Since the homology space $H_{2k}(Q^{2n-1}, \mathbb Z)$ of the
quadric $Q^{2n-1}$ of odd dimension is $\mathbb Z[X_w]$, any
subvariety of dimension $k$ has homology class $r[X_w]$ for some
$r \in \mathbb Z$. So  no Schubert variety in $Q^{2n-1}$ is Schur
rigid. If $X_w$ is a proper linear subspace $\PP^k$ $(k <n)$ of a
linear space $\PP^n$ contained in $G/P$, then $\CZ_k(G/P, r[X_w])$
contains $\CZ_k(\PP^n, r[\PP^k])$ so $X_w$ is not Schur rigid.

When $X_w$ is a maximal linear space in Hermitian symmetric space
$G/P$ of classical type or is a sub-Lagrangian Grassmannian
$L_{m-a}$ in the Lagrangian Grassmannian $L_m$, the Schur rigidity
of $X_w$ is proved in \cite{B}(Theorem 2, 5, 14, 15, 18 of
\cite{B}). Generalizing this result, we prove that, with the above
trivial exceptions, smooth Schubert varieties are Schur rigid.

\vs {\bf Main  Theorem.} {\it Let $G/P$ be a Hermitian symmetric
space other than  an odd dimensional quadric. Then any smooth
Schubert variety $X_w$ in $G/P$ is Schur rigid except when $X_w$
is a non-maximal linear space in $G/P$. Here, we consider $G/P$ as
a projective variety by the minimal equivariant embedding $G/P
\subset \PP(V)$. } \vs

After giving some preliminaries on Schubert varieties in Section
2,  we find a criterion that is used for determining  the Schubert
rigidity of a Schubert variety $X_w$ and one for the equality
$\CB_w=\CR_{w}$ in Section 3.

To prove the Schubert rigidity, we use the result of Goncharov on
the integral varieties of $F$-structures(\cite{G}). The fiber
$B_w$ of $\CB_w$ is an orbit of the semisimple part of  $P$, which
can be considered as a subgroup of $GL(k, \fm)$, where $\fm \simeq
T_0(G/P)$. In this case, Goncharov showed that if the cohomology
space $H^{1,1}(B_w) $ vanishes, then there is only one integral
varieties of $\CB_w$ passing through a fixed point and tangent to
a fixed tangent subspace. Thus $X_w$ is Schubert rigid if
$H^{1,1}(B_w)$ vanishes(Section \ref{Fstr}).

For the equality $B_w=R_w \subset Gr(k, \fm)$ of the fibers of
$\CB_w$ and $\CR_w$ at $o \in G/P$, we consider the decomposition
$ \wedge^k \fm=\bigoplus_{\ell(w)=k} \bI_w $ of the $k$-wedge
product of the tangent space $T_o(G/P) \simeq \fm$ as the direct
sum of irreducible  representation spaces of the semisimple part
of $P$ (\cite{K2}) and observe that $R_w$ is given by the
intersection $Gr(k, \fm) \cap \PP(\bI_w)$ and that $B_w$ is the
highest weight orbit in $\PP(\bI_w)$. Thus  if the complement of
tangent space $T_{[\fn_w]}B_w$ in $T_{[\fn_w]}Gr(k, \fm)$
intersects the tangent space $T_{[\wedge^k \fn_w]} \PP(\bI_w)$
trivially, then $B_w$ equals to $R_w$(Section \ref{S32}).

To compute the cohomology space $H^{1,1}(B_w)$ and the complement
of $T_{[\fn_w]}B_w$, we use the fact that a smooth Schubert
variety $X_w$ in a Hermitian symmetric space corresponds to a
subdiagram $\delta$ of the Dynkin diagram of $G$(\cite{BP},
\cite{LW}, Proposition \ref{smoothdynkin}).  Then we can apply the
theory of Lie algebra cohomology developed by Kostant (\cite{K1})
to compute them(Proposition \ref{SchubertD} and Proposition
\ref{equalityD}).

In Section 4, we verify these conditions and prove the main
Theorem. Furthermore, we prove the Schubert rigidity of smooth
Schubert varieties in general homogeneous space in the same way
under some assumptions(Proposition \ref{generalhomogeneous}).

Our representation theoretic method of proving the equality
$B_w=R_w$ is new and can be applied to the case of singular
Schubert varieties in general, because the assumption that $X_w$
is smooth is used only for computing the complement of $T_oB_w$ in
$T_oGr(k, \fm)$ in a uniform way.
 We expect that the rigidity of smooth Schubert varieties will
 serve as the building blocks for proving the rigidity of singular
 Schubert varieties in general.

\vs \noindent {\bf Acknowledgments.} I would like to thank the
referee for several suggestions which improve this paper. This
work was supported by the Post-Doctorial Fellowship Program of
Korea Science and Engineering Foundation(KOSEF).

\section{Differential systems}

\subsection{Schubert varieties} 

Let $\fg$ be a complex simple Lie algebra. Choose a Cartan
subalgebra $\fh$ of $\fg$. Let $\Delta$ be the set of all roots of
$\fg$ with respect to $\fh$. Fix a system $\CS=\{ \alpha_1,
\cdots, \alpha_{n}\}$ of simple roots of $\fg$. Let $\Delta^+$ be
the set of positive roots with respect to $\CS$.

For a subset $I$ of $\CS$, set $\Delta_I=\Delta \cap \,\mathbb Z
I$, and let $\fp_I$ be the parabolic subalgebra generated by $I$,
that is, $\fp_I=\fp_0 + \fm^*$, where  $\fp_0:=\fh + \sum_{\alpha
\in \Delta_I} \fg_{\alpha}$ is the reductive part and
$\fm^*:=\sum_{\alpha \in \Delta^+ \backslash \Delta_I}$ is the
nilpotent part. Then we have $\fg=\fp_I + \fm$, where
$\fm=\sum_{\alpha \in \Delta^- \backslash \Delta_I}$. The empty
set $I=\emptyset$ corresponds to the Borel subalgebra
$\fb=\fh+\sum_{\alpha \in \Delta^+} \fg_{\alpha}$. Set
$\Delta(\fm^*)=\Delta^+ \backslash \Delta_I$.

Let $W$ be the Weyl group of $\fg$.  For each $w \in W$, define a
subset $\Delta(w)$ of $\Delta^+$ by $\Delta(w):=w \Delta^- \cap
\Delta^+$. Then the number $|\Delta(w)|$ of elements in
$\Delta(w)$ is equal to the length $\ell(w)$ of $w$. The subset
$W^{I}$ of $W$ defined by  $W^I=\{ w \in W : \Delta(w) \subset
\Delta(\fm^*) \}$ is equal to $\{w \in W : w^{-1}(\Delta_I^+)
\subset \Delta^+ \}$.

Define a dual set of weights by requiring $<\lambda_i,
\alpha_j^{\vee}>=\delta_{i,j}$, where
$\alpha_j^{\vee}=\frac{2\alpha_j}{<\alpha_j, \alpha_j>}$. Set
$\rho=\frac{1}{2}\sum_{\alpha \in \Delta^+} \alpha=\sum_{i}
\lambda_i$. Put $\rho_I=\sum_{\alpha_i \in \CS \backslash I}
\lambda_i$.  Then $W$ (respectively, $W^I=W/W_I$)  is bijective to
the orbit of $\rho$ (respectively, $\rho_I$) under the action of
$W$ given by $\rho \rightarrow w^{-1} \rho$(\cite{K2}, \cite{BE}).


\vs

Let $G$ be a connected Lie group with Lie algebra $\fg$ and let
$P_I$, $B$ and $H$ be Lie subgroups of $G$ corresponding to
$\fp_I$, $\fb$ and $\fh$. If $P=P_I$, then $\Delta_I$, $W^I$,
$\rho_I$ are also denoted by $\Delta_P$, $W^P$, $\rho_P$,
respectively.

The Weyl group $W$ of $\fg$ is isomorphic to the quotient
$N_G(H)/H$ of the normalizer of $H$ by $H$. So we may consider $w
\in W$ as an element in $G$. In the same way, if $U$ is a real
form of $G$, we may take $w$ in $U$, too.  Let $o$ be the origin
of $G/P$. For each $w \in W^P$, considering $w$ as an element of
$G$, let $V_w$ be the $B$-orbit $V_w=Bw^{-1}\cdot o$. Then $G/P$
is decomposed as a disjoint union $\cup_{w \in W^P} V_w$ of
$B$-orbits(Section 6 of \cite{K2}).

For $w \in W$, let $\fn^*_w$ (respectively, $\fn_w$) be the
nilpotent Lie subalgebra spanned by the root vectors in
$\Delta(w)$ (respectively, $-\Delta(w)$), and let $N_w^*$ and
$N_w$ be the subgroups corresponding to $\fn_w^*$ and $\fn_w$.
Since $w \Delta(w^{-1})=-\Delta(w)$, $N_{w^{-1}}^*$ and $N_w$ are
conjugate by an element in $G$, that is,
$wN_{w^{-1}}^*w^{-1}=N_w$.

Write $B=N^* H$, where $N^*$ is the unipotent part of $B$. Let
$\kappa$ be the elements in $W$ of maximal length. Then $\kappa
\Delta^- = \Delta^+$.
 The isotropy group of $N^*$
at $w^{-1} \cdot o$ is $N^*_{w^{-1} \kappa}$ and $V_w$ is
isomorphic to the $N_w$-orbit $N_w \cdot o$ at the origin by the
action of an element $w \in G$, that is, $w V_w = w N_{w^{-1}}^*
w^{-1} \cdot o = N_w \cdot o$(Section 6 of \cite{K2}). The closure
$X_w$ of $V_w$ is an irreducible subvarieties of dimension
$\ell(w)$. We call $X_w$ a {\it Schubert variety of type} $w$.



\begin{prop} \label{invform} Let $G/P$ be a  Hermitian symmetric space.
Let $U$ be a compact real form of $G$. Write $\fp=\fp_0+\fm^*$ and
$\fg=\fp + \fm$, where $\fp_0$ is the reductive part and $\fm^*$
is the nilpotent part. Then

(1) $\bigwedge^k\fm$ is decomposed as the direct sum $\bigoplus_{w
\in W^P(k)} (\bigwedge \fm)^{w \rho-\rho}$ of irreducible
$\fp_0$-representation spaces $\bI_w:=(\bigwedge \fm)^{w
\rho-\rho}$, where $W^P(k)=\{ w \in W^P: \ell(w)=k\}$. The highest
 weight vector in $\bI_w$ is the decomposable $\ell(w)$-vector
$e_{-\Delta(w)}$, the wedge product of root vectors $x_{\alpha}$
of roots $\alpha \in -\Delta(w)$.

(2) For $w \in W^P$, define $\phi_w$ by the $U$-invariant
differential $(k,k)$-form which is given by
$(\sqrt{-1})^{k^2}\sum_i \zeta_i \wedge \overline{\zeta}_i$ at the
origin, where $\{\zeta_i\}$ be an orthonormal basis of $\bI_w^*$.
Then $\phi_w$ is closed and  positive and satisfies

\begin{eqnarray*}<[\phi_v], [X_w]>=\int_{X_w} \phi_v=0, \text{ for } v
\not=w.\end{eqnarray*}
\end{prop}

\begin{pf}
(1) Theorem 5.14 or Corollary 8.2 of \cite{K1}.

(2) Section 5.6, Corollary 5.4 and Corollary 6.15 of \cite{K2}. To
see it more directly, note that $w V_w=N_w \cdot o$ and that we
can take $w$ in $U$. Since $\phi_v$ is $U$-invariant, we have
$\int_{V_w} \phi_v = \int_{N_w \cdot o} \phi_v=0$.
\end{pf}

\begin{rk}
Schubert varieties $X_w$ and closed  differential forms $\phi_w$
such that $\phi_v|_{X_w}=0$ for $v \not=w$,  are defined on
general homogeneous space $G/P$(Corollary 6.15 of \cite{K2}). But,
in general, $\phi_w$ is not of the form as in the Proposition
\ref{invform} and is more complicated.

\end{rk}

\begin{exam} Let $Gr(m,n)$ be the Grassmannian of $m$-dimensional subspaces of
$\BC^n$.
 For $\ba \in
P(m,n)=\{\ba=(a_1, \cdots, a_m): n-m \geq a_1 \geq \cdots \geq a_m
\}$, the {\it Schubert variety} $\sigma_{\ba}$ {\it of type} $\ba$
is defined by the set
$$\{ E \in Gr(m,n)| \,\,dim(E \cap \BC^{n-m+i-a_i}) \geq i \}.$$
$\sigma_{\ba}$ is a subvariety of $Gr(m,n)$ of codimension
$|\ba|:=a_1 +\cdots + a_m$.

For $\ba \in P(m,n)$, define its {\it dual} $\ba^*$ by
$$\ba^*=(n-m-a_m, \cdots, n-m-a_1)$$
\noindent and  define its {\it conjugate} $\ba'=(a'_1, \cdots,
a'_{n-m})$  by $$\ba'_i=\sharp \{ j | \,\,a_j \geq i \} \text{ for
} 1 \leq i \leq n-m.$$

Let $E \in Gr(m,n)$ and $Q$ be the quotient $V/E$. The $k$-th
wedge product of $T^*_E(Gr(m,n))=E \otimes Q^*$ is decomposed as
the direct sum of irreducible $SL(E) \times SL(Q^*)$
 representation spaces as follows(Exercise 6.11 of \cite{FH}):
 $$ \wedge^{k}(E
\otimes Q^*) =\bigoplus_{|\ba|=k} \BS_{\ba}(E) \otimes
\BS_{\ba'}(Q^*),$$ \noindent where $\BS_{\ba}$ is the Schur
functor of type $\ba$. Locally, $\phi_{\ba}$ can be written as the
sum $(\sqrt{-1})^{k^2}\sum_i \zeta_i \wedge \bar{\zeta}_i$, where
$\{\zeta_i\}$ is an orthonormal basis of $\BS_{\ba}(E) \otimes
\BS_{\ba'}(Q^*)$.
\end{exam}

\vs We will explain how to construct the Schur module (or Weyl
module) $\BS_{\ba}(V)$ for a vector space $V$ in  Section 4, where
we will prove the Schur rigidity.

\subsection{Schur and Schubert differential systems}

\begin{defn} Let $M$ be a manifold and let $Gr(k, TM)$ be the
Grassmannian bundle of $k$-subspaces of the tangent bundle $TM$. A
subvariety $\CF$ of $Gr(k, TM)$ is called a {\it differential
system} on $M$. A subvariety $X$ of $M$ is said to be an {\it
integral variety} of the differential system $\CF$ if at each
smooth point $x \in X$, the tangent space $T_xX$ is an element of
the fiber $\CF_x$. We say that $\CF$ is {\it integrable} if at
each point $x \in M$ and $y \in \CF_x$, there is an integral
variety passing through $x$ and tangent to the subspace $W_y$ of
$T_xM$ corresponding to $y$.
\end{defn}

In \cite{W} and \cite{B}, they consider two differential systems
$\CB_w$ and $\CR_w$ to solve the rigidity problem of the Schubert
variety $X_w$. By taking different Borel subgroup $B \subset P$,
we get a family of Schubert varieties of type $w$. This family
induces the first differential system.

\begin{defn} Let $M=G/P$ be a Hermitian symmetric space.
 For each $w \in W^P$, the {\it Schubert differential
system} $\CB_{w}$ of type $w$ is the differential system with the
fiber consisting of the tangent spaces of all the Schubert
varieties $X_w$ of type $w$ passing through the given point. We
say that $X_w$ is {\it Schubert rigid} if Schubert varieties of
type $w$ are the only irreducible integral varieties of $\CB_{w}$.
\end{defn}

\begin{defn} Let $M=G/P$ be a Hermitian symmetric space.
For each $w \in W^P$, the {\it Schur differential system}
$\CR_{w}$ of type $w$ is the differential system with the fiber
defined by the intersection
$$\cap_{v \not=w, \ell(v)=\ell(w)} Z(\phi_{v}),$$
\noindent where $Z(\phi_{v})$ is the set of $\ell(w)$-subspace of
$T_xM$ on which $\phi_{v}$ vanishes. We say that $X_w$ is {\it
Schur rigid} if Schubert varieties  of type $w$ are the only
irreducible  integral varieties of $\CR_{w}$. \end{defn}

Proposition \ref{invform} provides a way to transform the problem
on the cycle space $\mathcal Z(G/P, r[X_w])$ to the problem on the
integral varieties of the differential system $\mathcal R_w$.
(Section 2.8.1 of \cite{B}, \cite{W}).

\begin{prop} Let $G/P$ be a Hermitian symmetric space. Take $w \in W^P$.
 For a subvariety $X$ of $G/P$, $X$ has the homology class
$[X]=r[X_{w}]$ for an integer $r$ if and only if $X$ is an
integral variety of $\CR_{w}$. Therefore, the cycle space
$\CZ(G/P, r[X_w])$ consists of formal sums of Schubert varieties
of type $w$ if and only if $X_w$ is Schur rigid.
\end{prop}


 For $w \in W^P$,
let $B_{w}$ (resp. $R_{w}$) be the fiber of $\CB_{w}$ (resp.
$\CR_{w}$) at $o \in G/P$.
 Note that $B_{w}$
is closed and $R_{w}$ is connected (Remark 2 and Remark 12 of
\cite{B}). Proposition \ref{invform} gives more refined structure
of $R_{w}$ and $B_{w}$.

\begin{prop} \label{hwo}  For $w \in W^P$,

 (1) $R_{w}$  is equal to the intersection $$Gr(k, \fm) \cap
\PP(\bI_{w}) \subset \PP(\wedge^{k}\fm)$$ and

 (2) $B_{w}$ is the  orbit of a highest weight vector in the
irreducible representation space $\PP(\bI_w)$ of the reductive
component $P_0$ of $P$.
\end{prop}

\begin{pf} (1) Consider a complex $k$-subspace $W$ of $\fm \simeq T_{x}(G/P)$
as a one-dimensional vector space $\wedge^k W$ in $\wedge^k \fm$
which is decomposed as $\bigoplus_{\ell(w)=k} \bI_w$. By
Proposition \ref{invform} (2),
 $\phi_{v}$ is given by the sum $(\sqrt{-1})^{k^2}\sum_i \xi_i
\wedge \overline{\xi}_i$ at $x$, where $\{\xi_i \}$ is an
orthonormal basis of $\bI_v^*$. If $\phi_{v}$ vanishes on a
complex  $k$-subspace $W$ of $T_{x}(G/P)\simeq \fm$, then so does
any $\xi_{i}$. Thus every vector in the vector space $\wedge^k W$
is annihilated by any dual vector in the space $\bI_v^*$, so it is
contained in the complement $\bigoplus_{\tau \not= v} \bI_{\tau}$
of $\bI_v$ in $\wedge^k \fm$. Therefore, we have
$$Z(\phi_v)=Gr(k, \fm) \cap \PP( \bigoplus_{\tau \not= v,
\ell(\tau)=k} \bI_{\tau})$$ \noindent and hence the intersection
$R_w=\cap_{v \not=w, \ell(v)=k} Z(\phi_v)$ is equal to $Gr(k, \fm)
\cap \PP(\bI_{w})$.

(2) Use Proposition \ref{invform} (1).
\end{pf}

\begin{rk} In \cite{B}, the author computes the integral elements of
the exterior differential ideal $\mathcal I_w$ which is generated
by the sections of the subbundle of $\wedge^k T^*(G/P)$ associated
to the component $\bI_w^*$ in $\wedge^k \fm^*$ and then computes
$R_w=\cap_{v \not=w, \ell(v)=\ell(w)} Z(\phi_v)$ by using the fact
that $Z(\phi_v)$ is equal to the set of all $k$-dimensional
integral element of $\mathcal I_v$(Lemma 1 in \cite{B}). Since
$Z(\phi_v)$ is given by the intersection $Gr(k, \fm) \cap \PP(
\bigoplus_{\tau \not= v, \ell(\tau)=k} \bI_{\tau})$, computing it
is much more difficult than computing just their intersection
$R_w=Gr(k, \fm) \cap \PP(\bI_{w})$. He uses it to classify
holomorphic vector bundles which are globally generated and
certain Chern classes of which are vanishing. In this paper we are
interested only in the rigidity problem so it is enough to
consider the intersection $R_w=Gr(k, \fm) \cap \PP(\bI_{w})$.
\end{rk}

\begin{prop} \label{equivalence}
Let $G/P$ be a Hermitian symmetric space. Take $w \in W^P$. Then
$X_w$ is Schur rigid if and only if $B_w$ is equal to $R_w$ and
$X_w$ is Schubert rigid.
\end{prop}

\begin{pf}
If $R_w$ is equal to $B_w$ and $X_w$ is Schubert rigid, then $X_w$
is Schur rigid by definition.  If $X_w$ is not Schubert rigid,
then $X_w$ is not Schur rigid because $B_w$ is a subvariety of
$R_w$.

Suppose that $B_w$ is not equal to $R_w$. Then we can find an
irreducible subvariety of $G/P$ which is an integral variety of
$\CR_w$ but which is not an integral variety of $\CB_w$ in the
same way as Bryant does in the case when $G/P$ is the Grassmannian
$Gr(m,n)$(Remark 4 and Proposition 2 of \cite{B}). We will present
the way to get such a subvariety for completeness.

For a subspace $A$ of $\fm \simeq T_x(G/P) $ with $[A] \in R_w
\backslash B_w$, consider the closure $X_A=cl(exp(A).x)$ of the
orbit of the unipotent group $exp(A)$ in $G/P$. Since $\fm$ is
abelian Lie algebra, any subspace $A$ of $\fm$ is again an abelian
Lie algebra and thus is a nilpotent subalgebra of $\fm$. So
$exp(A)$ is nothing but $Id + A$ if we consider it in the matrix
algebra.

Then $exp(A).x$ is Zariski  open dense in its closure $X_A$ so
$X_A$ is an irreducible algebraic variety of $G/P$. Since the form
$\phi_v$ vanishes on $exp(A).x$ for $v \not= w$, $X_A$ is an
integral variety of $\CR_w$. But its tangent space at $x$ is not
contained in $B_w$, so $X_w$ is not a Schubert variety of type
$w$. Hence $X_w$ is not Schur rigid.
\end{pf}

\subsection{Smooth Schubert varieties} Let $\CD(G)$ be the
Dynkin diagram of $G$.  Let $\gamma$ denote the simple root such
that the complement $\CS \backslash \{\gamma \}$ generates the
parabolic subgroup $P$. Then the marked diagram $(\CD(G), \gamma)$
represents the Hermitian symmetric space $G/P$.

To  a connected subgraph $\delta$ of $(\mathcal D(G), \gamma)$
with $\gamma$ a node we associate a smooth Schubert variety of
type, say, $w \in
 W^P$. Here, $w$ is the longest element in the subgroup of the Weyl group
 of $G$ generated by the reflections by
 the simple roots in $\delta$ and $\Delta(w^{-1})$ is equal to the set of all roots which
 are linear combinations of simple roots in $\delta$ with positive coefficients.

Conversely, we will show that in Hermitian symmetric space any
smooth Schubert variety corresponds to a subgraph of the Dynkin
diagram $\CD(G)$ of $G$. To do this we need the following
Proposition on the singular locus of Schubert varieties in
Hermitian symmetric spaces.

\begin{prop} Let $X_w$ be a Schubert variety in a Hermitian
symmetric space $G/P$.  Then the stabilizer $P_w$ in $G$ of $X_w$
is the parabolic subgroup generated by $\CS \cap w^{-1}(\Delta_P
\cup \Delta^{-})$. Furthermore, $X_w$ is smooth if and only if
$P_w$ acts on $X_w$ transitively.

\end{prop}

\begin{pf}
Section 2.6  and Proposition 3.3 and Proposition 4.4 of \cite{BP}.
Note that $W^P$ in \cite{BP} is the set of the inverse $w^{-1}$ of
$w$ in $W^P$ in our paper.
\end{pf}

\begin{prop} \label{smoothdynkin} Let $G/P$ a Hermitian symmetric
space and $(\mathcal{D}(G), \gamma)$ be the marked Dynkin diagram
correspond to $G/P$. Any smooth Schubert variety $X_w$ in $G/P$
corresponds to a connected subgraph of the Dynkin diagram
$\mathcal{D}(G)$ with $\gamma$ a node.
\end{prop}

\begin{pf}
Let $P_w=L_w U_w$ be the Levi decomposition of $P_w$, where $L_w$
is the Levi part and $U_w$ is the unipotent part. Then the
semisimple part of $L_w$ has the simple root system $\CS \cap
w^{-1}(\Delta_P \cup \Delta^{-})$ and $U_w$ is generated by
$U_{\alpha}$, for $\alpha \in \Delta^+ \backslash w^{-1}(\Delta_P
\cup \Delta^{-})$.

The isotropy of $N^*$ at  $w^{-1} \cdot o$ is $N^*_{w^{-1}
\kappa}$ and $\Delta(w^{-1} \kappa)=\Delta^{+} \cap w^{-1} \kappa
\Delta^- = \Delta^+ \cap w^{-1} \Delta^+$(Section 6 of \cite{K2}).
Thus $U_w$ and semisimple Lie subgroup of $L_w$ corresponding to
$\CS \cap w^{-1}(\Delta^+_P)$ act trivially on $w^{-1} \cdot o$.
So $X_w$ corresponds to the subgraph $\delta$ of $\CD(G)$ having
$\CS \cap w^{-1} \Delta^{-}$ as the set of nodes.
\end{pf}

For classical group $G$, we can verify Proposition
\ref{smoothdynkin} case-by-case by using the description of smooth
Schubert varieties in Section 5 of \cite{LW}.

\vs Now consider $G/P$ as the orbit of the highest weight vector
of the irreducible representation space $\PP(V^{\rho_P})$ with
highest weight $\rho_P$. This is the minimal equivariant embedding
of $G/P$ in the projective space. Then

(1)  a subgraph of $\CD(G)$ of type $A_k$ with $\gamma$ an
extremal node corresponds to a linear space $\PP^{k}$ in $G/P$. A
subgraph of $\CD(G)$ of type $A_{\ell}$ with $\gamma$ a node in
the middle corresponds to a subgrassmannian. A Schubert variety of
this type exists when $G$ is of type $A_n, D_n, E_6, E_7$.

(2) A subgraph of $\CD(G)$ of type $D_{\ell}$ with $\gamma$ an
extremal node corresponds to a quadric or an isotropic
Grassmannian in $G/P$. A Schubert variety of this type exists when
$G$ is of type $D_n, E_6, E_7$.

(3) A subgraph of $\CD(G)$ of type $C_{\ell}$ with
$\gamma=\alpha_{\ell}$ (in the usual notation of simple roots)
corresponds to Lagrangian Grassmannian $L_n$ in $G/P$. A Schubert
variety of this type exists when $G/P$ is Lagrangian Grassmannian
$L_m$.

These are all smooth Schubert varieties in $G/P$(Proposition
\ref{smoothdynkin}).

\section{Criterions for the rigidity}

\subsection{Schubert rigidity} \label{Fstr}
Schubert differential systems are $F$-structures for various $F$'s
and the Schubert rigidity can be proved by showing the vanishing
of a certain cohomology space.

\begin{defn}
 Let $F $  be a submanifold of $ Gr(k, V)$
with a transitive action of a subgroup of $GL(V)$, where  $\dim
V=n$. A fiber bundle $\CF \subset Gr(k, TM)$ on a manifold $M$ of
dimension $n$ is said to be an {\it $F$-structure}   if  at each
point $x \in M$ there is a linear isomorphism $\varphi(x): V \ra
T_xM$ such that the induced map $\varphi(x)^k : Gr(k, V) \ra Gr(k,
T_xM)$ sends $F$ to $\CF_x$.
\end{defn}

One can get the information on the integrability and on the set of
integral varieties of an $F$-structure from certain cohomology
spaces $H^{k,\ell}(F)$ depending only on the embedding $F\subset
Gr(k,V)$(Chapter 1 of \cite{G}).

The first cohomology $H^{1,1}(F)$ is defined by
$$H^{1,1}(F)=Ker(
\partial :W_f^* \otimes T_fF \ra \wedge^2 W_f^* \otimes
(V/W_f)),$$ \noindent  where $W_f$ stands for the $k$-subspace of
$V$ represented by $f \in F$ and for $p: W_f \ra T_fF$, $\partial
p :\wedge^2 W_f \ra V/W_f$ is defined
 by $\partial p(V_1,
V_2)=p(V_1)(V_2)-p(V_2)(V_1)$, considering $T_fF$ as a subspace of
$W_f^* \otimes (V/W_f)$. Define $H^{j,1}(F)$ inductively so that
the vanishing of $H^{1,1}(F)$ imply the vanishing of $H^{j,1}(F)$
for all $j$ and that for a given integrable $F$-structure, if
$H^{n,1}(F)=0$, then the family of all integral varieties passing
through a fixed point and tangent to a fixed tangent subspace has
dimension $\sum_{j \leq n-1}\dim H^{j,1}(F)$(Chapter 1 of
\cite{G}). In short,

\begin{prop} \label{HF} Let $\CF$ be an integrable $F$-structure.
If $H^{1,1}(F)=0$, then for a fixed $x \in M$ and $y \in \CF_x$,
there exists a unique integral variety of $\CF$ passing through
$x$ and tangent to $W_y$.
\end{prop}

 Now we apply Proposition \ref{HF} to the Schubert differential system
$\CB_w$ on the Hermitian symmetric space $G/P$ as an
$B_w$-structure and . Since $P_0$ acts on the fiber $B_w$
transitively, we have the decomposition $\fp_0=\fm_w + \fl_w +
\fm_w^*$, where $\fl_w+\fm_w^*$ is the Lie algebra of the isotropy
group. Then the tangent space $T_{[\fn_w]}B_w \subset
T_{[\fn_w]}Gr(k,\fm)$  is isomorphic to $\fm_w \subset \fn_w^*
\otimes \fm/\fn_w$.

\begin{prop}
Let $G/P$ be a Hermitian symmetric space. Fix $w \in W^P$. Let
$B_w$ be the fiber of the Schubert differential system $\CB_w$. If
$H^{1,1}(B_w)=Ker(\partial : \fn_w^* \otimes \fm_w \rightarrow
\wedge^2 \fn_w^* \otimes \fm/\fn_w)$ is zero, then $X_w$ is
Schubert rigid.
\end{prop}

\subsection{The equality $B_w=R_w$} \label{S32}

 To find when $B_w$ is equal to $R_w$, we will
compare their tangent spaces.

Note that under the embedding $Gr(k,\fm) \subset \PP(\wedge^k
\fm)$,  a tangent vector $\varphi$ in $T_{[\fn_w]}Gr(k,\fm)\simeq
\fn_w^* \otimes \fm/\fn_w$ can be considered as a tangent vector
$\varphi^k$ in $T_{[\wedge^k \fn_w]}\PP(\wedge^k \fm) \simeq
\wedge^k \fn_w^* \otimes \wedge^k\fm/\wedge^k \fn_w$ in the
following way. Take a basis $\{v_1, \cdots, v_k \}$ of $\fn_w$.
For a linear map $\varphi:\fn_w \ra \fm/\fn_w$, the map
$\varphi^k: \wedge^k \fn_w \ra \wedge^k \fm/\wedge^k\fn_w$ is
given by
$$\varphi^k(v_1 \wedge \cdots \wedge v_k)=\sum_i v_1 \wedge \cdots
\wedge \varphi(v_i) \wedge \cdots \wedge v_k \mod \wedge^k \fn_w.
$$

\begin{prop} \label{comparetangent}
Suppose that for the highest weight vector $\varphi$ of every
irreducible $\fl_w$-representation space in the complement of
$\fm_w$ in $\fn_w^* \otimes \fm/\fn_w$,
$$\varphi^k(\wedge^k \fn_w) \not \subset \bI_w/\wedge^k \fn_w. $$
\noindent   Then $B_w$ is equal to $R_w$.

\end{prop}

\begin{pf}
Put
$$T_w=\{ \varphi:\fn_{w} \ra \fm/\fn_w | \,\,
\varphi^k(\wedge^k \fn_w) \subset \bI_w /\wedge^k \fn_w \}.$$
  Since the space  $T_w$ is equal to
the intersection $T_{[\fn_{w}]} Gr(k, \fm) \cap \,T_{[\wedge^k
\fn_w]}\PP(\bI_w)$ and each space is invariant by the action of
$\fl_w$, $T_w$ is an $\fl_w$-representation space. The tangent
space $T_{[\fn_w]} B_w=\fm_w$ is contained in $T_w$. By the
hypothesis, the complement of $\fm_w$ in $T_{[\fn_w]}Gr(k, \fm)$
intersects $T_{[\wedge^k \fn_w]}\PP(\bI_w)$ trivially. Thus the
tangent space $T_{[\fn_w]} B_w =\fm_w$ is equal to $T_w$ in $
\fn_w^* \otimes \fm/\fn_w$. Then since the tangent space $\fm_w$
of $B_w$ is contained in the tangent space $T_{[\wedge^k \fn_{w}]}
R_{w}$ of $R_{w}$ and the latter is contained in
 $T_w$, $B_w$ and $R_w$ have the same tangent space at $[ \fn_w]$.
Hence  $B_w$ is equal to $R_{w}$.
\end{pf}

\subsection{Lie algebra cohomology}

The complement of $\fm_w$ in $\fn_w^* \otimes \fm/\fn_w$, together
with $H^{1,1}(B_w)$, can be computed using Lie algebra cohomology
as follows.

\begin{prop}    \label{liecohomology} For $w \in W^P$,
define the action $\rho$ of $\fn_{w}$ on the space $\fm_{w} +
\fm/\fn_{w}$ by

\begin{displaymath}
\rho(A)(X) = \left\{ \begin{array}{ll}
 X(A) \in \fm/\fn_{w} & \text{ for } X \in \fm_{w} \subset
\fn_{w}^* \otimes
(\fm/\fn_{w}), \\
0 & \text{ for } X \in \fm/\fn_{w}.
\end{array} \right.
\end{displaymath}

\noindent Then the Lie algebra cohomology $H^1(\fn_{w}, \fm_{w} +
\fm/\fn_{w})$ associated to the representation of the nilpotent
Lie algebra $\fn_{w}$ on $\fm_{w} + \fm/\fn_{w}$ is the direct sum
of $H^{1,1}(B_w)=Ker (\partial : \fn_{w}^* \otimes \fm_{w} \ra
\wedge^2 \fn_{w}^* \otimes (\fm/\fn_{w}))$ and the complement of
$\fm_w$ in $\fn_w^* \otimes \fm/\fn_w$.
\end{prop}

\begin{pf}
The complex defining the Lie algebra cohomology $H(\fn_{w},\fm_{w}
+ \fm/\fn_{w})$
$$ \fm_{w} + \fm/\fn_{w}\ra
 \fn_{w}^* \otimes (\fm_{w} + \fm/\fn_{w})
\ra  \wedge^2 \fn_{w}^* \otimes (\fm_{w} + \fm/\fn_{w} ) \ra \quad
\cdots
$$ is the direct sum of the following subcomplexes:

\begin{eqnarray*}
(1) \quad \qquad 0 \quad \ra \quad  \fn_{w}^* \otimes \fm_{w}
\quad \quad &\ra& \wedge^2 \fn_{w}^*
\otimes (\fm/\fn_{w}) \ra \quad \cdots \\
(2)\quad \qquad \, \,\fm_{w} \ra\quad  \fn_{w}^* \otimes
(\fm/\fn_{w}) &\ra& \qquad 0 \qquad \ra \quad \cdots
\end{eqnarray*}

Since $[\fn_{w}, \fn_{w}]=0$, the boundary map of the first
complex ($1$) is the same as the map $\partial$.
\end{pf}

In general, computing the Lie algebra cohomology $H^1(\fn_w, \fm_w
+ \fm/\fn_w)$ is not easy. However, if there is a semisimple Lie
algebra $\fq_0$ which contains $\fn_w$, such that

(1) $\fn_w$ is the maximal nilpotent ideal of a parabolic
subalgebra of $\fq_0$

 (2) the representation of $\fn_w$ on $\fm_w + \fm/\fn_w$ is the
 restriction of a representation of $\fq_0$ on $\fm_w +
 \fm/\fn_w$,

 \noindent then it can be calculated by the work of Kostant(\cite{K1}).
When the Schubert variety $X_w$ is associated with a subgraph of
the Dynkin diagram $\CD(G)$ of $G$, we can find such a Lie algebra
$\fq_0$ in a natural way, which we will be treated in the
following section.

\subsection{Subdiagrams of the marked Dynkin diagram $(\CD(G), \gamma)$}

\label{subdiagram}
 Fix  a connected subgraph
$\delta$ of $(\mathcal D(G), \gamma)$ with $\gamma$ a node. Let
$X_w$ be a smooth Schubert variety corresponding to $\delta$. In
this case, the Schubert rigidity of $X_w$ can be checked by the
same method as in \cite{CH}, where they deal with the cases when
the Schubert varieties are linear spaces after  $G/P$ being
embedded in the projective space minimally. In the following we
collect the facts that are needed to deal with the cases studied
in this paper. \vs

From now on, to avoid negative signs in the computation, we will
adopt the convention that nilpotent part of Borel or parabolic
subalgebras are generated by   root spaces of negative roots. Then
the tangent space of the homogeneous space $G/P$ is generated by
root spaces of positive roots.

For $\alpha \in \CS$, let $\sigma_{\alpha}$ denote the reflection
with respect to $\alpha$, which is an element of the Weyl group of
$\fg$. For a nilpotent Lie algebra $\fn$ of $\fg$, let
$\Delta(\fn)$ denote the set of all roots $\beta$ with
$\fg_{\alpha} \subset \fn$. Then $\fn$ can be written as
$\fn=\sum_{\alpha \in \Delta(\fn)} \fg_{\alpha}$.

For $\alpha_i \in \mathcal S$, define $n_{\alpha_i}:\Delta
\rightarrow \mathbb Z$ by $n_{\alpha_i}(\alpha)=n_i$ for
$\alpha=\sum_{j} n_j \alpha_j$.
 For a subset $\mathcal S_1=\{\alpha_{i_1},
\cdots, \alpha_{i_r}\}$ of $\CS$, define the map $n_{\mathcal
S_1}$ by
\begin{eqnarray*}
n_{\mathcal S_1}: \,\, \Delta^+ &\longrightarrow& \BZ^{|I|} \\
        \alpha=\sum n_i \alpha_i &\longmapsto& ( n_{i_1}, \cdots,
        n_{i_r}).
\end{eqnarray*}

 For a connected subgraph $\delta$
of $\CD(G)$, define the neighborhood $N(\delta)$ of $\delta$ by
the set of simple roots that are not in $\delta$ but that are
connected to $\delta$ by an edge.

\vs

Consider  the map $n_{\CS_1}$ for $\CS_1=N(\delta)$.  For each
$\bi=(n_{i_1}, \cdots, n_{i_r})$ in the image
$n_{N(\delta)}(\Delta^+) \subset \mathbb Z^{|N(\delta)|}$, let
$\lambda_{\bi}$ denote the root in $\{ \alpha \in \Delta ^+:
n_{N(\delta)}(\alpha)=\bi\}$ of maximal height. Let $D$ be the set
of all such $\lambda_{\bi}$.

\begin{prop} \label{lac}
Let $\delta$ be a connected subgraph of $(\CD(G), \gamma)$ with
$\gamma$ a node. Let $X_w$ be a schubert variety of type $w$
corresponding to $\delta$. Define a reductive Lie algebra $\fq_0$
by $\fq_0=\fh+ \sum_{\alpha \in \Delta_{\delta}} \fg_{\alpha}$,
where $\Delta_{\delta}=\Delta_{\CS \backslash N(\delta)}$ and
define $D$ as above.
 Then $H^1(\fn_w, \fm_w+
\fm/\fn_w)$ is the direct sum of irreducible
$\fq_0$-representation spaces $H^1(\fn_w, \fm_w+
\fm/\fn_w)_{\lambda}$, $\lambda \in D$. The highest weight vector
of $H^1(\fn_w, \fm_w+ \fm/\fn_w)_{\lambda}$ is
$x_{\gamma}^*\otimes x_{\sigma_{\gamma}(\lambda)}$.
\end{prop}

\begin{pf} Set $\fn=\sum_{\alpha \in \Delta^+ \backslash
\Delta_{\delta}} \fg_{\alpha}$. Then $\fq_0+ \fn^*$ is the
parabolic subalgebra generated by $\CS \backslash N(\delta)$ and
we have $\fg=\fn+ \fq_0+ \fn^*$. The nilpotent Lie algebras $\fn,
\fm, \fn_w, \fm_w$ are generated by the root spaces of roots in
\begin{eqnarray*}
\Delta(\fn)&=&\{ \alpha \in \Delta^+: n_{\beta}(\alpha) >0 \text{
for some } \beta \in N(\delta) \}
\\
\Delta(\fm)&=&\{ \alpha \in \Delta^+: n_{\gamma}(\alpha)=1 \}
\\
\Delta(\fm_w)&=& \{ \alpha \in \Delta^+ : n_{\gamma}(\alpha)=0,
\,\, n_{\beta}(\alpha) > 0 \text{ for some } \beta \in N(\delta)
\} \\
\Delta(\fn_w) &=& \{ \alpha \in \Delta^+ : n_{\gamma}(\alpha)=1,
\,\, n_{\beta}(\alpha)=0 \}.
\end{eqnarray*}
So  we have $\fn=\fm_w+\fm/\fn_w$ and $\fq_0=\fn_w + \fl_w
+\fn_w^*$. Furthermore, the restriction of the adjoint action of
$\fq_0$ on $\fn$ to $\fn_w$ is equal to the action of $\fn_w$ on
$\fm_w + \fm/\fn_w$ defined in Proposition \ref{liecohomology}.


For each $\bi=(n_{i_1}, \cdots, n_{i_r})$ in the image
$n_{N(\delta)}(\Delta^+) \subset \mathbb Z^{|N(\delta)|}$, let
$\fn_{\bi}$ denote the nilpotent Lie subalgebra spanned by the
root vectors $x_{\alpha}$  of root $\alpha$ with
$n_{N(\delta)}(\alpha)=\bi $.  Then $\fn$ is decomposed as
$\sum_{\bi} \fn_{\bi}$ and each $\fn_{\bi}$ is an irreducible
$\fq_0$-representation space with highest weight
$\lambda_{\bi}$(Section 9.9 of \cite{BE}).
 Since the action by $\fn_w$
does not change the coefficient $n_{\beta}(\alpha)$ for $\beta \in
N(\delta)$, $H^1(\fn_w, \fn)$ is the direct sum of  $H^1(\fn_w,
\fn_{\bi})$'s.
 By Theorem 5.14 of \cite{K1}, each $H^1(\fn_w,
\fn_{\bi})$ is an irreducible representation space of $\fl_w$ with
the highest weight vector  $x_{\gamma}^* \otimes
x_{\sigma_{\gamma}(\lambda_{\bi})}$.
\end{pf}

\begin{prop} \label{SchubertD} With the same notations as in Proposition \ref{lac}, put
\begin{eqnarray*}
 D'&=& \{ \lambda \in D : n_{\gamma}(\sigma_{\gamma}(\lambda))=0
 \} \\
 &=& \{ \lambda \in D : \sum_{\beta \in N(\gamma)}
 n_{\beta}(\lambda)=n_{\gamma}(\lambda) \}.
\end{eqnarray*}

 \noindent Then $H^{1,1}(B_w)$ is the direct summand of $H^1(\fn_w,
\fm_w+\fm/\fn_w)$
 consisting of $\fl_w$-modules with the highest weight vector
 $x_{\gamma}^* \otimes x_{\sigma_{\gamma}(\lambda)}$ for $\lambda \in
 D'$.
Therefore, $X_w$ is Schubert rigid  if $D'=\emptyset$.
\end{prop}

\begin{pf}
It suffices to show that $n_{\gamma}(\sigma_{\gamma}(\lambda))=0$
if and only if $\sum_{\beta \in N(\gamma)}
n_{\beta}(\lambda)-n_{\gamma}(\lambda)=0$.  By definition,
$\sigma_{\gamma}(\alpha)=\alpha -<\alpha, \gamma^{\vee}> \gamma$.
For a simple root $\alpha$,
\begin{eqnarray*}
<\alpha, \gamma^{\vee}>= \begin{cases}\,\, 2 \text{ if } \alpha=\gamma \\
-1 \text{ if } \alpha \in N(\gamma) \\
\,\, 0 \text{ otherwise }\end{cases}
\end{eqnarray*}
so we have $n_{\gamma}(\sigma_{\gamma}(\lambda))=\sum_{\beta \in
N(\gamma)} n_{\beta}(\lambda)-n_{\gamma}(\lambda)$. This gives
the desired equality.
\end{pf}

\begin{prop}  \label{equalityD} With the same notations as in Proposition \ref{lac}, put
\begin{eqnarray*}
 D''&=& \{ \lambda \in D : n_{\gamma}(\sigma_{\gamma}(\lambda))=1
 \} \\
 &=& \{ \lambda \in D : \sum_{\beta \in N(\gamma)}
 n_{\beta}(\lambda)=n_{\gamma}(\lambda)+1 \}.
\end{eqnarray*}
Take a basis $\{ v_1, \cdots, v_k\}$ of $\fn_w$ consisting of
weight vectors and with $v_1=x_{\gamma}$. If
$x_{\sigma_{\gamma}(\lambda)} \wedge v_2 \wedge \cdots \wedge v_k$
is not contained in $\bI_w$ for any $\lambda \in D''$, then $B_w$
is equal to $R_{w}$.

\end{prop}

\begin{pf}
By Proposition \ref{comparetangent} and Proposition
\ref{liecohomology}, to show that $\fm_w$ is equal to $T_w$, it
suffices to check for these generators $x_{\gamma}^* \otimes
x_{\sigma_{\gamma} (\lambda)}$ with $\lambda \in D''$ that the
image of the map $\varphi_{\lambda}^k$ induced by
$\varphi_{\lambda}=x_{\gamma}^* \otimes
x_{\sigma_{\gamma}(\lambda)}$ is not contained in $\bI_w$. Since
we take $v_1=x_{\gamma}$,  $\varphi_{\lambda}^k(v_1 \wedge \cdots
\wedge v_k)=x_{\sigma_{\gamma}(\lambda)} \wedge v_2 \cdots \wedge
v_k$. \end{pf}

\vs The Schubert differential system $\CB_w$ and the Schubert
rigidity can be defined in general homogeneous space $G/P$ in the
same way as in the case of Hermitian symmetric spaces. Suppose
that $P$ is a maximal parabolic corresponding to a simple root
$\gamma$. To a subdiagram $\delta$ of $\CD(G)$ with $\gamma$ a
node associate  a Schubert variety $X_w$ in $G/P$ in the same way.
Assume that the marked diagram $(\delta, \gamma)$ represents a
Hermitian symmetric space. Then the tangent space $\fn_w$ of $X_w$
satisfies $[\fn_w, \fn_w]=0$. So we can apply the same method as
above to check the Schubert rigidity of such a Schubert variety.
Note that $\fm$ is not irreducible: $\fm=\sum \fm_i$, where
$\fm_i=\sum_{n_{\gamma}(\alpha)=i}\fg_{\alpha}$. However, $\fn_w$
is contained in $\fm_1$ and the complex defining Lie algebra
cohomology $H(\fn_{w},\fm_{w} + \fm/\fn_{w})$ is decomposed as
follows:

\begin{eqnarray*}
(1) \quad \qquad 0 \quad \ra \quad  \fn_{w}^* \otimes \fm_{w}
\quad \quad &\ra& \wedge^2 \fn_{w}^*
\otimes (\fm_1/\fn_{w}) \quad \ra \quad \cdots \\
(2)\,\,\quad \quad \, \,\fm_{w} \ra\quad  \fn_{w}^* \otimes
(\fm_1/\fn_{w}) &\ra& \qquad  \wedge^2\fn_w^* \otimes \fm_2 \qquad
\ra
\quad \cdots \\
(3) \qquad \quad  \, \,\fm_1 \ra\quad  \fn_{w}^* \otimes \fm_2
\qquad
&\ra& \qquad  \wedge^2\fn_w^* \otimes \fm_3 \qquad \ra \quad \cdots \\
...
\end{eqnarray*}
Thus $X_w$ is Schubert rigid  if $D'=\emptyset$ as in Proposition
\ref{SchubertD}.

\section{Schubert rigidity and Schur rigidity}

\subsection{Schubert rigidity} \label{kGP}

 Let $\delta_0$ be the subgraph of $\CD(G)$ of type
$A_{\ell}$ with two ends, $\gamma$ and the branch point(the node
connected to three other points). Let $k_{G/P}$ be the length of
$\delta_0$. Then for any $j \leq k_{G/P}$, $H^{2j}(G/P, \mathbb Z)
\simeq \mathbb Z$ and thus the homology class of any subvariety of
dimension $j \leq k_{G/P}$ is a multiple of the homology class of
the Schubert variety $\PP^j, j \leq k_{G/P}$ in $G/P$. So these
Schubert varieties are not Schur rigid.

\begin{prop} \label{SchubertHSS} Let $G/P$ be a Hermitian symmetric space
which is not an odd dimensional quadric.
 Let $(\delta, \gamma)$ be a subgraph of $(\CD(G), \gamma)$
 of type $(A_k, \alpha_1)$ with length $> k_{G/P}$ and let
 $\delta_{max}$ be the unique maximal one among subgraphs of $(\CD(G), \gamma)$
 of type $(A_n, \alpha_1)$ containing
 $\delta$.

 Let $X_w$ be a Schubert variety corresponding to
 $\delta$.
 Then  integral varieties of
 $\CB_w$ are subvarieties of a Schubert variety of type
 corresponding to $\delta_{max}$. Therefore $X_w$ is Schubert rigid
 if $X_w$ is a maximal linear space.
\end{prop}

\begin{pf} See \cite{CH}. \end{pf}

\begin{lem} \label{empty} Let $(\CD(G), \gamma)$ be a marked
Dynkin diagram and $(\delta, \gamma)$ be a subgraph.  For $\bi \in
n_{N(\delta)}(\Delta^+)$, let $\lambda_{\bi}$ be the root of
maximal height  among $\alpha \in \Delta^+$ with $
n_{N(\delta)}(\alpha)=\bi$. If $(\delta, \gamma)$ is not of type
$(A_k, \alpha_1)$, then
$n_{\gamma}(\sigma_{\gamma}(\lambda_{\bi}))>0 $
 for all $\bi$.
\end{lem}

\begin{pf} From the proof of Proposition \ref{SchubertD} it follows that
 $n_{\gamma}(\sigma_{\gamma}(\lambda_{\bi}))\geq \sum_{\beta \in
N(\gamma)} n_{\beta}(\lambda_{\bi})-n_{\gamma}(\lambda_{\bi})$. If
$(\delta, \gamma)$ is not of type $(A_k, \alpha_1)$, then for
every $\bi \in n_{N(\delta)}(\Delta^+)$, we have $\sum_{\beta \in
N(\gamma)} n_{\beta}(\lambda_{\bi}) >n_{\gamma}(\lambda_{\bi})$
because $\lambda_{\bi}$ has maximal height among $\alpha \in
\Delta^+$ with $ n_{N(\delta)}(\alpha)=\bi$.
\end{pf}

\begin{prop} \label{smoothSchubert} Let $G/P$
be a  Hermitian symmetric space. Let $(\delta, \gamma)$ be a
subgraph of $(\CD(G), \gamma)$ whic is not
 of type $(A_k, \alpha_1)$. Then the corresponding Schubert variety $X_w$
of $G/P $ is Schubert rigid.
\end{prop}

\begin{pf} Use Proposition
\ref{SchubertD} and Lemma \ref{empty}.
\end{pf}


Lemma \ref{empty} holds for general homogeneous space. By the
remark at the end of Section \ref{subdiagram}, we have

\begin{prop} \label{generalhomogeneous}
Let $G/P$ be a homogeneous space with $P \subset G$ a maximal
parabolic subgroup corresponding to $\gamma$. Let $\delta$ be a
subgraph of $\CD(G)$ with $\gamma$ a node. Supposed that the
marked diagram $(\delta, \gamma)$ represents a Hermitian symmetric
space but that it is not of type $(A_k, \alpha_1)$. Then the
corresponding Schubert variety $X_w$ in $G/P$ is Schubert rigid.
\end{prop}

\begin{exam} Schubert varieties as in Proposition
\ref{generalhomogeneous} arise when we  consider the
 desingularization of  singular Schubert
varieties in Hermitian symmetric space.

Let $G/P$ be an even dimensional quadric $Q^{2n} \subset \PP(V)$
and $X_w$ be a singular Schubert variety of codimension 1. Then
$X_w$ may be thought as the locus of the family of  $\PP^1$'s in
$Q^{2n}$ passing through a fixed point $p \in Q$, or,
equivalently, of the family of  $\PP(E)$'s  for isotropic
2-subspaces $E$ of $V$ containing a fixed isotropic 1-subspace
$W$.

For an isotropic 2-subspace $E$ of $V$, $\PP(E)$ is a smooth
Schubert variety in $Q^{2n}$, say, of type $v$.  The family of
isotropic 2-subspaces of $V$ is a homogeneous space, say, $G/Q$.
Then the subvariety of $G/Q$ consisting of isotropic 2-subspaces
containing a fixed isotropic 1-subspace $W$ is a smooth Schubert
variety $X_u$ of $G/Q$ of the form as in Proposition
\ref{generalhomogeneous}. Thus we can write $X_w=\cup_{W \subset
E} \PP(E)$ as the locus $\cup_{t \in X_u} (X_v)_t$ of the family
$\{ (X_v)_t \}$ of smooth Schubert varieties of type $v$ in $G/P$
parameterized by a smooth Schubert variety $X_u$ of type $u$ in
$G/Q$.

In general, a singular Schubert variety $X_w$ in $Q^{2n}$ of
codimension $r <n-1$ can be expressed as the locus of the family
of smooth Schubert varieties in $Q^{2n}$ parameterized by a smooth
Schubert varieties in $G/Q$ of the form as in Proposition
\ref{generalhomogeneous}. \vs

$X_w$, a singular Schubert variety of codimension $r$ in
$G/P=Q^{2n}$

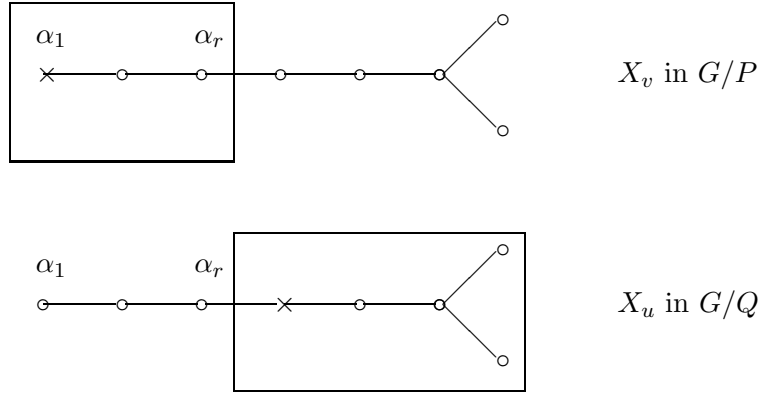
\vs
\begin{picture} (200, 50)

\put(0, 15){$\alpha_1$} \put(60,15){$\alpha_r$}

\put(0,0) {$\times$}  \put(3,3){\line(1,0){27}}
\put(30,0){$\circ$} \put(34, 3){\line(1,0){27}}

\put(60,0) {$\circ$}  \put(64,3){\line(1,0){27}}

\put(90,0) {$\circ$}  \put(94,3){\line(1,0){27}}

\put(120,0) {$\circ$}  \put(124,3){\line(1,0){27}}

\put(150,0) {$\circ$}  \put(154,3){\line(1,1){20}}

\put(150,0) {$\circ$}  \put(154,3){\line(1,-1){20}}

\put (174, 21) {$\circ$}
 \put (174, -21) {$\circ$}



\put(220,0){$X_v$ in $G/P$}

\put(-10, 30){\line(1,0){85}} \put(-10, -30){\line(1,0){85}}

\put(-10, 30){\line(0,-1){60}} \put(75, 30){\line(0,-1){60}}

\end{picture}

\vs \vs \vs

\begin{picture} (200, 50)

\put(0, 15){$\alpha_1$} \put(60,15){$\alpha_r$}

\put(90,0) {$\times$}  \put(3,3){\line(1,0){27}}
\put(30,0){$\circ$} \put(34, 3){\line(1,0){27}}

\put(60,0) {$\circ$}  \put(64,3){\line(1,0){27}}

\put(0,0) {$\circ$}  \put(94,3){\line(1,0){27}}

\put(120,0) {$\circ$}  \put(124,3){\line(1,0){27}}

\put(150,0) {$\circ$}  \put(154,3){\line(1,1){20}}

\put(150,0) {$\circ$}  \put(154,3){\line(1,-1){20}}

\put (174, 21) {$\circ$}
 \put (174, -21) {$\circ$}



\put(220,0){$X_u$ in $G/Q$}

\put(75, 30){\line(1,0){110}} \put(75, -30){\line(1,0){110}}

\put(75, 30){\line(0,-1){60}} \put(185, 30){\line(0,-1){60}}

\end{picture}

\vs \vs \vs \vs \vs

A desingularization of of $X_w$ can be obtained by considering the
following double fibration:
\begin{eqnarray*}
Gr_{iso}(r+1, V)=G/Q \stackrel{\mu}{\longleftarrow}
Flag_{iso}(1,r+1, V) \stackrel{\nu}{\longrightarrow} Q^{2n}=G/P
\end{eqnarray*}
Fix an isotropic $r$-subspace $W$ in $V$. Then the union of
$\PP^r$'s in $Q^{2n}$ containing $\PP(W)$ is a Schubert variety
$X_w$ of type $w$ in $Q^{2n}=G/P$ of codimension $r$. The
subvariety
 consisting of isotropic $(r+1)$-subspaces of $V$ containing
$W$  is a Schubert variety $X_u$ of type $u$ in $Gr_{iso}(r+1,
V)=G/Q$. The preimage $\mu^{-1}(X_u) \subset Flag_{iso}(1,r+1, V)$
is smooth and  maps to $X_w$ birationally by $\nu$.
\end{exam}

\subsection{The equality $B_w=R_w$: $Gr(m,n)$-case}

First, we review how to construct the Schur module (or  Weyl
module) $\BS_{\ba}(V)$ for a vector space $V$ of dimension $n$ and
for a partition $\ba=(a_1, \cdots, a_n)$. For references, see
Chapter 6 and Chapter 15 of \cite{FH} and Chapter 8 of \cite{F}.

Let $\ba'$ be the conjugate of $\ba$ and $k=|\ba|$. Put the
factors of the $k$-th tensor product $V^{\otimes k}$ in one-to-one
correspondence with the squares of the Young diagram of $\ba$.
Then $\BS_{\ba}(V)$ is the image of this composite map:
$$ \otimes_i(\wedge^{a'_i}V) \ra \otimes_i (\otimes^{a'_i}V) \ra
V^{\otimes k} \ra \otimes_j (\otimes^{a_j} V) \ra
\otimes_j(Sym^{a_j}(V)),$$ the first map being the tensor product
of the obvious inclusion, the second grouping the factors of
$V^{\otimes k}$ according to the columns of the Young diagram, the
third grouping the factors according to the rows of the Young
diagram, and the forth being the obvious quotient map(Exercise
6.14 of \cite{FH}).

Let $\{v_1, \cdots, v_n \}$ be a basis of $V$. For a semi-stable
tableau $T$ of shape $\ba$ (a numbering of the Young diagram of
$\ba$ in such a way that the number is strictly increasing in
column and is nondecreasing in row), let $e_T \in
\otimes_i(\wedge^{a'_i}V)$ be a tensor product of wedge products
of basis element for $V$, the $i$-th factor in $\wedge^{a'_i}V$
being the wedge product (in order) of those basis vectors whose
indices occur in the $i$-th column of $T$. Then the set of the
images of these elements $e_T$ for semi-stable tableau $T$, is a
basis of $\BS_{\ba}(V)$(Theorem 1 of Section 8.1 of \cite{F}).
Furthermore, under the action of $SL(V)$, each $e_T$ is a weight
vector of weight $\bb=(b_1, \cdots, b_m)$, where $b_i$ is the
number of times the integer $i$ occurs in $T$(Section 8.2 of
\cite{F}). The highest weight vector corresponds to the
semi-standard tableau $T$ with all the elements in $i$-th row
equal to $i$. e.g.
 in $\BS_{\ba}(V)$ for $\ba=(4,4,4,4)$,
\vs
\begin{tabular}{|c|c|c|c|} \hline 1 &1 & 1& 1 \\ \hline 2& 2& 2& 2 \\ \hline 3
&3&3& 3 \\ \hline
\end{tabular} corresponds to the weight vector of highest weight  $(4,4,4,4)$ \vs
and \vs
\begin{tabular}{|c|c|c|c|} \hline 1 &2 & 2&3 \\ \hline 2& 3& 3&4 \\ \hline 3
&4&5&5 \\ \hline
\end{tabular} corresponds to the weight vector of weight
$(1,3,4,2,2)$.

\vs

\vs The dimension of the weight space with weight $\bb$ in the
representation space $\BS_{\ba}(\BC^n)$ is the number  of ways one
can fill the Young diagram of $\ba$ with $b_1$ $1$'s, $b_2$ $2$'s,
$\cdots$ $b_n$ $n$'s, in such a way that the entries in row are
nondecreasing and those in each column are strictly increasing.
e.g. For $\ba=(4,4,4)$, $\bb=(4,3,2,2,1)$ is a weight in the
representation space $\BS_{\ba}(V)$ but $\bb=(5,4,3)$ is not.

\begin{prop} \label{smoothSchur}
For  $\ba=(p^q)$ with $(p,q) \not=(1,1)$, the Schubert
differential system $\CB_{\ba^*}$ is equal to the Schur
differential system $\CR_{\ba}$.
\end{prop}

\vs  We remark that the case when $p=1$ or $q=1$  was proved in
Lemma 4 and Lemma 5 of \cite{B}.

\vs

\begin{pf}  In this case
$\sigma_{\ba^*}$ is $\{ E \in Gr(m,n): \BC^{m-q} \subset E \subset
\BC^{m+p} \} \simeq Gr(q, p+q)$. For  simplicity, assume that $1 <
p < n-m$ and $1 < q < m$(other cases can be proved in a similar
way). Let $\delta$ be the subdiagram consisting of the simple
roots $\alpha_i$ for $ m-q +1 \leq i \leq m+p-1$. Then
$N(\delta)=\{ \alpha_{m-q}, \alpha_{m+p} \}$ and $D=\{
\lambda_{1,0}, \lambda_{0,1}, \lambda_{1,1} \}=D''$, where
$\lambda_{1,0}=\sum_{1 \leq i \leq m+p-1} \alpha_i$,
$\lambda_{0,1}=\sum_{m-q+1 \leq i \leq n} \alpha_i$, and
$\lambda_{1,1} =\sum_{1 \leq i \leq n} \alpha_i$. For all $\lambda
\in D''$, we have $\sigma_{\gamma}(\lambda)=\lambda$.

Recall that $\bI_{\ba}$ is given by $\BS_{\ba}(E^*) \otimes
\BS_{\ba'}(Q)$. Choose a basis $\{e_1, \cdots, e_m \}$
(respectively, $\{ q_1, \cdots, q_{n-m} \}$) of $E$ (respectively,
of $Q$) such that the vector $e_i^* \otimes q_j$ is a root vector
of root $\sum_{m-i+1 \leq s \leq m+j-1} \alpha_s$ in $\fm=E^*
\otimes Q$. The basis of $\BS_{\ba}(E^*) \otimes \BS_{\ba'}(Q)$
consisting of weight vectors is
\begin{eqnarray*} \{ e_T^* \otimes q_S  \,\,: &&T
\text{ is a semi-standard tableau of shape } \ba \text{ and }
\\ &&S \text{ is a semi-standard tableau of shape } \ba'
 \},
\end{eqnarray*}
\noindent where $e_T$ and $q_S$ are defined at the beginning of
this section.


The vector space $\fn_{\ba^*}$ is spanned by the vectors $e_i^*
\otimes q_j$, $1 \leq i \leq q, 1 \leq j \leq p$ and $\wedge^k
\fn_{\ba^*}=\wedge_{i,j}(e_i^* \otimes q_j)$ is a scalar multiple
of $e_T^* \otimes q_S$ for  the semi-standard tableau \vs
$T=$\begin{tabular}{|c|c|c|c|}
\hline 1&1& $\cdots$ &1 \\
\hline 2 & 2& $\cdots$ &2 \\
\hline $\vdots$ & $\vdots$ &  &$\vdots$ \\
\hline $q$ & $q$ & $\cdots$ & $q$ \\ \hline
\end{tabular}
and $S^t=$\begin{tabular}{|c|c|c|c|}
\hline 1&2& $\cdots$ &$p$ \\
\hline 1 & 2& $\cdots$ &$p$ \\
\hline $\vdots$ & $\vdots$ &  &$\vdots$ \\
\hline $1$ & $2$ & $\cdots$ & $p$ \\ \hline
\end{tabular}
\vs \noindent Here, we denote by $S^t$ the transpose of $S$, by
noting that the Young diagram of $\ba'$ is the transpose of the
Young diagram of $\ba$. This will help to understand the element
$e_T^* \otimes q_S$ in the tensor product $\BS_{\ba}(E^*) \otimes
\BS_{\ba'}(Q)$. Also, note that $e_T^* \otimes q_S$ is not a
decomposable $k$-vector in $\wedge^k(E^* \otimes Q)$ for general
$T$ and $S$ but that, for $T$ and $S$ as above, $e_T ^* \otimes
q_S$ is a highest weight vector and it is a decomposable
$k$-vector $\wedge_{i,j} (e_i^* \otimes q_j)$.

 The weight vector $x_{\lambda_{1,0}}$ of weight
$\lambda_{1,0}$ is $e_{m}^* \otimes q_{p}$. If we replace $e_1^*
\otimes q_1$ by $e_{m}^* \otimes q_{p}$, then
$\varphi_{1,0}^k:=(e_{m}^* \otimes q_{p}) \wedge (e_1^* \otimes
q_2) \wedge \cdots \wedge (e_q^* \otimes q_p)$ is a weight vector
of weight $(p-1, p,\cdots, p,1)\oplus(q-1, q, \cdots, q, q+1)$.
But $(q-1, q, \cdots, q,q+1)$ is not a weight of $\BS_{\ba'}(Q)$.
Thus $\varphi_{1,0}^k$ is not contained in
$\bI_{\ba}=\BS_{\ba}(E^*) \otimes \BS_{\ba'}(Q)$. Similarly,
$\varphi_{0,1}^k=x_{\lambda_{0,1}} \wedge (e_1^* \otimes q_2)
\wedge \cdots \wedge (e_q^* \otimes q_p)$ is not contained in
$\bI_{\ba}$.

For $\lambda_{1,1}$, consider $\varphi_{1,1}^k=x_{\lambda_{1,1}}
\wedge (e_1^* \otimes q_2) \wedge \cdots \wedge (e_q^* \otimes
q_p)$, where $x_{\lambda_{1,1}}=e_m^* \otimes q_{n-m}$ is the
weight vector of weight $\lambda_{1,1}$. Since
$ad(x_{-\alpha_{m+p}-\cdots-\alpha_n})$ does not change any vector
in $\varphi^k_{1,1}$ except $x_{\lambda_{1,1}}$, we get
$$ad(x_{-\alpha_{m+p}-\cdots-\alpha_n})\varphi^k_{1,1}=x_{\lambda_{1,0}}
\wedge (e_1^* \otimes q_2) \wedge \cdots \wedge (e_q^* \otimes
q_p) \not\in \bI_{\ba}.$$ Since
$ad(x_{-\alpha_{m+p}-\cdots-\alpha_n})$ preserves each irreducible
$\fp_0$-representation in $\wedge^k \fm$, $\varphi^k_{1,1}$ is not
contained in $\bI_{\ba}$. By Proposition \ref{equalityD},
$B_{\ba^*}$ is equal to $R_{\ba}$.
\end{pf}

\subsection{The equality $B_w=R_w$: Hermitian symmetric space
$\not=Gr(m,n)$ case}

To prove the equality $B_w=R_w$ for the smooth Schubert variety
$X_w$ in Hermitian symmetric space $G/P \not= Gr(m,n)$, we will
consider a natural irreducible representation space $\bI_{\tau}$
in $\wedge^k \fm$ other than $\bI_w$ so that  for any $\lambda \in
D''$, $\varphi^k_{\lambda}=x_{\sigma_{\gamma}(\lambda)} \wedge v_2
\wedge \cdots \wedge v_k$ has a nonzero $\bI_{\tau}$-component,
i.e., if we write $\varphi^k_{\lambda}=\sum_{\sigma}
\varphi_{\sigma}$, where $\varphi_{\sigma} \in \bI_{\sigma}$
according to the decomposition $\wedge^k \fm=\sum_{\sigma}
\bI_{\sigma}$, then $\varphi_{\tau} \not=0$. So it will not be
contained in $\bI_w$. This will prove the equality $B_w=R_w$ by
Proposition \ref{equalityD}.

Note that if $w=\sigma_{\beta_1} \cdots \sigma_{\beta_k}$,
$\beta_i \in \CS$, is a reduced expression of $w$, then
$\sigma_{\beta_1} \cdots \sigma_{\beta_{i-1}}(\beta_i), i=1,
\cdots, k$ are all distinct and $\Delta(w)$ consists of such
roots.

\begin{prop} \label{SchurHSS} Let $G/P$ be a Hermitian symmetric space
$\not=Gr(m,n)$. Let $\delta$ be a subgraph of $\CD(G)$ of type
$A_{\ell}, \ell > k_{G/P}$ with $\gamma$ a  node and $X_w$ be the
corresponding Schubert variety. Then $\CB_w$ is equal to $\CR_w$.
\end{prop}

Here, $k_{G/P}$ is the length of the chain connecting $\gamma$ and
the branch point as in Section \ref{kGP}.

\vs  We remark that the case of the maximal linear spaces in the
quadric $Q^{2n}$ and the case of  linear spaces of dimension $\geq
3$ in the isotropic Grassmannian $N^+_{2n}$ was proved in Lemma 12
and Lemma 13 of \cite{B}.

\vs

\begin{pf} First, we will prove for the case of quadric $Q^{2n}$ and then
 show that other cases can be reduced to the quadric case.

{\bf Special case. } Assume that $G/P$ is $Q^{2(n-1)}$. Let
$\CS=\{ \alpha_1, \cdots, \alpha_n \}$ be the set of simple roots
indexed in the following way.

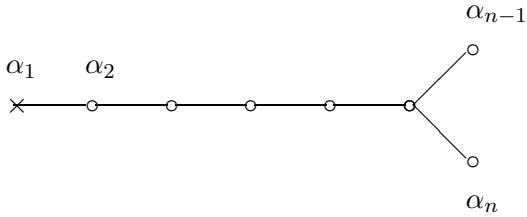
\vs
\begin{picture} (200, 50)

\put(0, 15){$\alpha_1$} \put(30,15){$\alpha_2$}

\put(0,0) {$\times$}  \put(3,3){\line(1,0){27}}
\put(30,0){$\circ$} \put(34, 3){\line(1,0){27}}

\put(60,0) {$\circ$}  \put(64,3){\line(1,0){27}}

\put(90,0) {$\circ$}  \put(94,3){\line(1,0){27}}

\put(120,0) {$\circ$}  \put(124,3){\line(1,0){27}}

\put(150,0) {$\circ$}  \put(154,3){\line(1,1){20}}

\put(150,0) {$\circ$}  \put(154,3){\line(1,-1){20}}

\put (174, 21) {$\circ$}
 \put (174, -21) {$\circ$}

\put(174, 36){$\alpha_{n-1}$} \put (174, -36){$\alpha_n$}

\end{picture}

\vs \vs \vs \vs \vs

 A maximal linear space $\PP^{n-1}$ in $Q^{2(n-1)}$ corresponds
to the element $w=\sigma_1 \cdots \sigma_{n-2} \sigma_{n-1}$ or
the element $\tau=\sigma_1 \cdots \sigma_{n-2}\sigma_n$ in $W^P$.
Here, $\sigma_i$ stands for the reflection with respect to the
simple root $\alpha_i$.

Consider the Schubert differential system $\CB_w$ and Schur
differential system $\CR_w$ for the Schubert variety $X_w$. Then
\begin{eqnarray*}
\Delta(w) &=& \{\alpha_1, \,\, \alpha_1 + \alpha_2,\,\, \cdots,
\,\, \alpha_1 + \alpha_2 + \cdots + \alpha_{n-2} + \alpha_{n-1} \}
\\
D &=& \{ \lambda_1:=\alpha_1 + 2 \alpha_2 + \cdots + 2
\alpha_{n-2} + \alpha_{n-1} + \alpha_n \}. \end{eqnarray*} and
thus $s_{\sigma_{\gamma}(\lambda_1)}=x_{\lambda_1}$. We will  show
that $$\varphi^k_{\lambda_1}=x_{\lambda_1} \wedge x_{\alpha_1+
\alpha_2} \wedge \cdots \wedge x_{\alpha_1 + \alpha_2 + \cdots +
\alpha_{n-2} + \alpha_{n-1}}$$ is contained in
$\bI_{\tau}=(\wedge^{n-1}\fm)^{\rho-\tau \rho}$.

The lowest weight vector of $\bI_{\tau}$ is given by $$
e_{\Delta(\tau)}=x_{\alpha_1} \wedge x_{\alpha_1+ \alpha_2} \wedge
\cdots \wedge x_{\alpha_1+ \alpha_2 + \cdots +\alpha_{n-2}} \wedge
x_{\alpha_1+ \alpha_2 + \cdots \alpha_{n-2} + \alpha_n}.$$ The
action by the root group $U_{\alpha_2 + \cdots + \alpha_{n-2}
+\alpha_{n-1}}$ does not change any vector in the multivector
$e_{\Delta(\tau)}$ except the first and the last one.  Applying
$U_{\alpha_2 + \cdots + \alpha_{n-2} +\alpha_{n-1}}$ to
$e_{\Delta(\tau)}$, we get an element
$$ x_{\alpha_1 + \alpha_2 + \cdots \alpha_{n-2} + \alpha_{n-1} }
\wedge x_{\alpha_1 + \alpha_2} \wedge \cdots \wedge x_{\alpha_1+
\alpha_2 + \cdots + \alpha_{n-2}} \wedge x_{\lambda_1}$$ in
$\bI_{\tau}$, which is a scalar multiple of
$\varphi^k_{\lambda_1}$. Hence $\varphi^k_{\lambda_1}$ is not
contained in $\bI_w$ and, by Proposition \ref{equalityD}, $\CB_w$
is equal to $\CR_w$.

\vs {\bf General case. } Let $w \in W^P$ correspond to the
subgraph $\delta$ of type $A_{k}$ with $k=\ell(w)$. Since $\gamma$
is an end, we may choose the simple root system $\CS=\{
\alpha_1=\gamma, \,\, \alpha_2, \cdots, \alpha_n \}$ such that $w$
is written by $\sigma_1 \cdots \sigma_k$. Put $k_0=k_{G/P}$. Let
$\delta_0$ be the subgraph of $\CD(G)$ of type $A_{k_0}$ with two
ends, $\gamma$ and the branch point (=the node connected to three
other nodes).

 Let $\epsilon$ be the subgraph of  $\CD(G)$ of type
$D_{k_0+1}$ which contains $\delta_0$. We may assume that the
simple roots are indexed in such a way that $\{\alpha_1, \cdots,
\alpha_{k_0}, \alpha_{k_0+1}, \alpha_{k+1} \}$ is the set of
simple roots in the subgraph $\epsilon$.

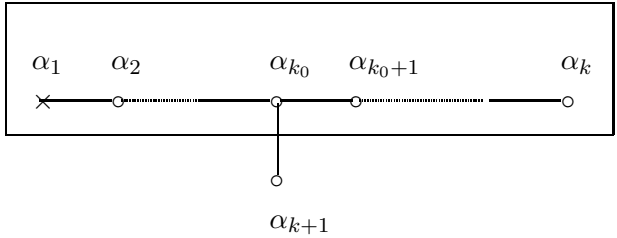
\vs
\begin{picture} (200, 50)
\put(-10, 40){\line(1,0){230}} \put(-10, 40){\line(0,-1){50}}
\put(-10, -10){\line(1,0){230}} \put(220, 40){\line(0,-1){50}}

\put(0, 15){$\alpha_1$} \put(30,15){$\alpha_2$}

\put(0,0) {$\times$}  \put(3,3){\line(1,0){27}}
\put(30,0){$\circ$}

\multiput(34, 3) (1,0){47} {\line(1,0){0.2}}

\put(90, 15){$\alpha_{k_0}$} \put(90, -45){$\alpha_{k+1}$}

 \put(90,0) {$\circ$}

\put(120, 15){$\alpha_{k_0+1}$}
 \put(120, 0) {$\circ$}
\multiput(124, 3) (1,0){47} {\line(1,0){0.2}}

\put(90, 3) {\line (-1,0){27}} \put(94,3){\line(1,0){27}}

 \put(93, 2){\line(0,-1){27}} \put(90, -30){$\circ$}

\put(200, 3){\line(-1,0){27}} \put(200, 0){$\circ$}
\put(200,15){$\alpha_k$}
\end{picture}

\vs \vs \vs \vs \vs

Put $\tau=\sigma_1 \sigma_2 \cdots  \sigma_{k-1} \sigma_{k+1}
\not=w \in W^P$. Since $\Delta(\tau)$ can be obtained from
$\Delta(w)$ by replacing $ \theta_{\delta}:=\alpha_1 + \cdots +
\alpha_k $ by $\theta_{\tau}:=\alpha_1 + \cdots + \alpha_{k_0} +
\alpha_{k+1}$, the lowest weight vector of $\bI_{\tau}$ is given
by
$$e_{\Delta(\tau)}=x_{\alpha_1} \wedge \cdots  \wedge
 \hat{x}_{\theta_{\delta}} \wedge
x_{\theta_{\tau}}. $$

\noindent Here, we use the notation $\hat{x}_{\theta_{\delta}}$ to
emphasize that $e_{\Delta(\tau)}$ does not have
$x_{\theta_{\delta}}$ in its decomposition, compared to
$e_{\Delta(w)}$ which can be considered as a base multivector. We
will use this kind of notation to express various multivectors.

 Write  $\theta_{\delta}=\alpha_1 +\beta_1$. Then
$\beta_1=\alpha_2 + \cdots + \alpha_k$ is a root in $\delta$, in
particular, in $\fl_w$. The action of the root group $U_{\beta_1}$
does not change any vector in the multivector $e_{\Delta(\tau)}$
except $x_{\alpha_1}$ and $x_{\theta_{\tau}}$. Applying
$U_{\beta_1}$ to $e_{\Delta(\tau)}$, we get an element
$$\hat{x}_{\alpha_1} \wedge \cdots \wedge x_{\theta_{\delta}}
\wedge x_{\theta_{\tau}+ \beta_1} \in \bI_{\tau}$$

 We will show that for $\lambda \in D''$,
$\varphi^k_{\lambda}=\hat{x}_{\alpha_1} \wedge \cdots \wedge
x_{\theta_{\delta}} \wedge x_{\lambda}$ has a nonzero
$\bI_{\tau}$-component. We will prove it by showing that, after a
successive adjoint action by the root vectors in $\fp_0$ to
$\varphi^k_{\lambda}$, we get  $\hat{x}_{\alpha_1} \wedge \cdots
\wedge x_{\theta_{\delta}} \wedge x_{\theta_{\tau}+ \beta_1} \in
\bI_{\tau}$.

First, assume that $\lambda$ is the minimal one in $D''$. Since
the coefficient in $\alpha_{k+1}$ of any $\lambda \in D''$ is not
zero, the $\lambda_{\bi}$ for $\bi=(1)$ (when $|N(\delta)|=1$) or
$\bi=(1,0)$ (when $|N(\delta)|=2$) is the minimal one in $D''$.
Then $\lambda-(\theta_{\tau} + \beta_1):=\beta_2$ is a root in
$\fp_0$.

If $\beta_2$ is not a root in $\fl_w$, then for any $\alpha \in
\Delta(w)$, $\alpha -\beta_2$ is not a root. If $\beta_2$ is a
root in $\fl_w$, then for some $\alpha \in \Delta(w)$, $\alpha -
\beta_2$ may be a root but this $\alpha-\beta_2$ is another root
in $\Delta(w)$. Since any $\lambda \in D''$ has coefficient $2$ in
$\alpha_2$, the coefficient in $\alpha_2$ of $\beta_2$ is zero and
thus  $\alpha_1 + \beta_2$ is not a root. So $\alpha-\beta_2$
cannot be $\alpha_1$ for any $\alpha \in \Delta(w)$. Thus, in any
case, $\hat{x}_{\alpha_1} \wedge \cdots \wedge [x_{-\beta_2},
x_{\alpha}] \wedge \cdots \wedge x_{\lambda} $ is zero for all
$\alpha \in \Delta(w)$. So we get
$$ad(x_{-\beta_2}) \varphi^k_{\lambda}= \pm \,\, \hat{x}_{\alpha_1}
\wedge \cdots \wedge x_{\theta_{\delta}} \wedge x_{\theta_{\tau} +
\beta_1} \in \bI_{\tau}.$$

\noindent  Since $ad(x_{-\beta_2})$ preserves irreducible
$\fp_0$-representation spaces in $\wedge^k \fm$,
$\varphi^k_{\lambda}$ has a nonzero $\bI_{\tau}$-component.

Other $\lambda$'s in $D''$ than the minimal one $\lambda_{min}$
can be obtained from the minimal one by adding a root
$\beta_{\lambda}$ in $\fp_0$ with nonzero coefficient in a simple
root in $N(\delta)$, successively. Then $\alpha-{\beta_{\lambda}}$
is not a root for any $\alpha \in \Delta(w)$. Thus we get
$$ad(x_{-\beta{\lambda}})\varphi^k_{\lambda} = \hat{x}_{\alpha_1} \wedge
\cdots \wedge x_{\theta_{\delta}} \wedge x_{\lambda_{min}} \not\in
\bI_w$$ and hence $\varphi^k_{\lambda} \not\in \bI_w$.
\end{pf}

\vs  Now we consider  the equality $B_w=R_w$ when the Schubert
variety $X_w$ corresponds to a subgraph $\delta$ of $\CD(G)$ of
type $D_{\ell}$ or $C_{\ell}$. The proof will be almost the same
but the difference is that in this case, the highest weight
$\theta_{\delta}$ is not obtained from $\alpha_1=\gamma$ by just
adding one root $\beta_1$ but by adding two roots $\beta_1$ and
$\beta_2$ successively. So we need to apply the adjoint actions
more times to get an element in $\bI_{\tau}\not= \bI_w$.

\begin{prop} \label{SchurHSSD} Let $G/P$ be a Hermitian symmetric space
$\not=Gr(m,n)$ associated to the simple root $\gamma$. Let
$\delta$ be a subgraph of $\CD(G)$ of type $D_{\ell}$ with
$\gamma$ an extremal node or of type $C_{\ell}$ with
$\gamma=\alpha_{\ell}$ and $X_w$ be the corresponding Schubert
variety. Then $\CB_w$ is equal to $\CR_w$.
\end{prop}

\vs  We remark that the case of sub-Lagrangian Grasmmanians in
Lagrangian Grassmannian $L_m$ was proved in Lemma 14 of \cite{B}.

\vs

\begin{pf} Let $\{\alpha_1, \cdots, \alpha_{\ell} \}$ be the simple
roots in $\delta$ with $\alpha_1=\gamma$. Then $w$ is the longest
element in the subgroup of the Weyl group generated by the
reflections $\sigma_i=\sigma_{\alpha_i}$, $i=1, \cdots, \ell$.
Write $w=\sigma_{i_1} \cdots \sigma_{i_k}$, where $k=\ell(w)$. Let
$\alpha_{\ell+1}$ be the simple root in $N(\delta)$. Put
$\tau=\sigma_{i_1} \cdots \sigma_{i_{k-1}} \sigma_{\ell+1}$. Then
$\ell(\tau)=\ell(w)$ and
$$\Delta(\tau)=\Delta(w) -\{\theta_{\delta}\} \cup \{
\theta_{\tau} \},$$ \noindent where $\theta_{\delta}$ is the
maximal root in the root system of $\delta$ and $\theta_{\tau}$ is
the maximal root in the root system of the subgraph $\epsilon$ of
type $A_{\ell}$ containing $\alpha_1$ and $\alpha_{\ell+1}$. After
rearranging the indexes of the simple roots, we may assume that
$\{\alpha_1, \cdots, \alpha_{\ell-1}, \alpha_{\ell+1} \}$ is the
simple root system of $\epsilon$. Then $\theta_{\tau}=\alpha_1+
\cdots + \alpha_{\ell-1} + \alpha_{\ell+1}$.

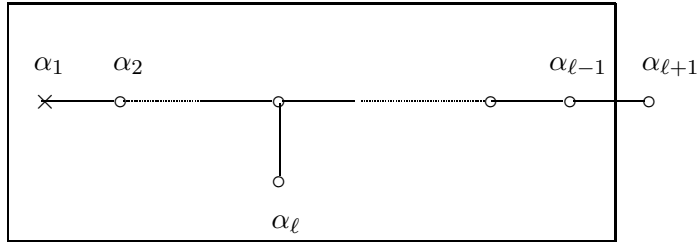
\vs
\begin{picture} (200, 50)
\put(-10, 40){\line(1,0){230}} \put(-10, 40){\line(0,-1){90}}
\put(-10, -50){\line(1,0){230}} \put(220, 40){\line(0,-1){90}}


\put(0, 15){$\alpha_1$} \put(30,15){$\alpha_2$}

\put(0,0) {$\times$}  \put(3,3){\line(1,0){27}}
\put(30,0){$\circ$}

\multiput(34, 3) (1,0){47} {\line(1,0){0.2}}


\put(90, -45){$\alpha_{\ell}$}

 \put(90,0) {$\circ$}

\multiput(124, 3) (1,0){47} {\line(1,0){0.2}}

\put(90, 3) {\line (-1,0){27}} \put(94,3){\line(1,0){27}}

 \put(93, 2){\line(0,-1){27}} \put(90, -30){$\circ$}

\put(170, 0){$\circ$}

 \put(200, 3){\line(-1,0){27}}

 \put(200,0){$\circ$} \put(195, 15){$\alpha_{\ell-1}$}

\put(204, 3){\line(1,0){27}} \put(230, 0){$\circ$} \put(230,
15){$\alpha_{\ell+1}$}
\end{picture}

\vs \vs \vs \vs \vs

We will show that for $\lambda \in D''$,
$\varphi^k_{\lambda}=\hat{x}_{\alpha_1} \wedge \cdots \wedge
x_{\theta_{\delta}} \wedge x_{\lambda}$ has a nonzero
$\bI_{\tau}$-component. By the same argument as in the proof of
Proposition \ref{SchurHSS}, it suffices to show this for the
minimal $\lambda_1$ in $D''$(In this case, since $|N(\delta)|=1$,
$\lambda_1$ is the minimal one).

 Write $\theta_{\delta}$ as the sum $ \alpha_1+ \beta_1+ \beta_2$
where $\beta_1:=\alpha_2+ \cdots + \alpha_{\ell-2}$ and
$\beta_2:=\alpha_2 + \cdots +\alpha_{\ell-1} + \alpha_{\ell}$ are
roots in $\delta$, in particular, in $\fl_w$.
 Let $\lambda=\lambda_1 $ be the minimal one in $ D''$.
Then one can check that $\lambda-\beta_1$ is not a root.

Write $\lambda=\theta_{\tau} + \beta_3 + \beta_4$. Here, $\beta_3$
is $\lambda - \theta_{\tau}$, if $\lambda-\theta_{\tau}$ is a
root, and $\beta_3$ is  the maximal root among the roots with each
coefficient $\leq$ each coefficient in $\lambda-\theta_{\tau}$,
otherwise. Then $\beta_4=0$ in the former case and
$\beta_4=\lambda-\theta_{\tau} -\beta_3$ in the latter case.

 Note that $\beta_3$ has
coefficient $1$ in $\alpha_{\ell}$, and $\beta_4$ is  a root in
$\fp_0$ with the coefficient $1$ in $\alpha_{\ell}$ and zero
coefficient in $\alpha_2$ when $\lambda-\theta_{\tau}$ is not a
root. So neither $\theta_{\tau}-\beta_2+ \beta_3$ nor
$\theta_{\tau}-\beta_2+\beta_4$ is a root(For,  the coefficient of
$\theta_{\tau}-\beta_2+ \beta_3$ and
$\theta_{\tau}-\beta_2+\beta_4$ in $\alpha_{\ell}$ is zero and any
root with the zero coefficient  in $\alpha_{\ell}$ has all the
coefficient $\leq 1$).

Since $\lambda-\beta_1$ is not a root,
 we get $$ad(x_{-\beta_1})(\varphi^k_{\lambda})= \pm
\,\,x_{\alpha_1} \wedge \cdots \wedge \hat{x}_{\alpha_1+\beta_1}
\wedge \cdots \wedge x_{\theta_{\delta}} \wedge x_{\lambda}$$
\noindent  so that
\begin{eqnarray*}&&ad(x_{-\beta_2})ad(x_{-\beta_1})(\varphi^k_{\lambda}) \\ &=& \pm
\,\,x_{\alpha_1}  \wedge \cdots \wedge \hat{x}_{\theta_{\delta}}
\wedge x_{\lambda} \pm \,\, x_{\alpha_1} \wedge \cdots \wedge
\hat{x}_{\alpha_1+\beta_1} \wedge \cdots \wedge
x_{\theta_{\delta}} \wedge [x_{-\beta_2}, x_{\lambda}].
\end{eqnarray*}
Since  $\theta_{\delta}$ is the maximal roots in $\delta$ and
$\lambda-\beta_2-\beta_3=\theta_{\tau}-\beta_2+\beta_4$ is not a
root, we get
\begin{eqnarray*} ad(x_{-\beta_3})ad(x_{-\beta_2})ad(x_{-\beta_1})
(\varphi^k_{\lambda})  =\pm \,\, x_{\alpha_1} \wedge \cdots \wedge
\hat{x}_{\theta_{\delta}} \wedge x_{\theta_{\tau} +\beta_4}+ {\bf
v},
\end{eqnarray*}
where $ {\bf v}$ is $\pm \,\, x_{\alpha_1} \wedge \cdots \wedge
\hat{x}_{\alpha_1+\beta_1+\beta_3} \wedge \cdots \wedge
x_{\theta_{\delta}} \wedge [x_{-\beta_2}, x_{\lambda }]$ if
$\alpha_1+\beta_1+\beta_3$  is a root in $\delta$, and is zero,
otherwise.
Since $\lambda-\beta_2-\beta_4=\theta_{\tau}-\beta_2+\beta_3$ is
not a root, we get
\begin{eqnarray*}ad(x_{-\beta_4})ad(x_{-\beta_3})ad(x_{-\beta_2})ad(x_{-\beta_1})(\varphi^k_{\lambda})
= \pm \,\,x_{\alpha_1}  \wedge \cdots \wedge
\hat{x}_{\theta_{\delta}} \wedge x_{\theta_{\tau}} + {\bf v'}
\end{eqnarray*}
\noindent where ${\bf v'}$ is a decomposable $k$-vector.

In any case, since $e_{\Delta(\tau)}=x_{\alpha_1} \wedge \cdots
\wedge \hat{x}_{\theta_{\delta}} \wedge x_{\theta_{\tau}}$ is a
$k$-vector in $\bI_{\tau}$ and $x_{-\beta_4}$, $x_{-\beta_3},
x_{-\beta_2}$ and $x_{-\beta_1}$ are root vectors of $\fp_0$,
$\varphi^k_{\lambda}$ has a nonzero $\bI_{\tau}$-component. Hence
$\varphi^k_{\lambda} \not\in \bI_w$ for any $\lambda \in D''$ and
thus $B_w=R_w$ by Proposition \ref{equalityD}.

\vs The case when $\delta$ is of type $(C_{\ell}, \alpha_{\ell})$
can be proved in a similar way by taking $\tau=\sigma_1 \cdots
\sigma_{i_{k-1}} \sigma_{\ell+1}$ after re-index the simple roots
in the following way.

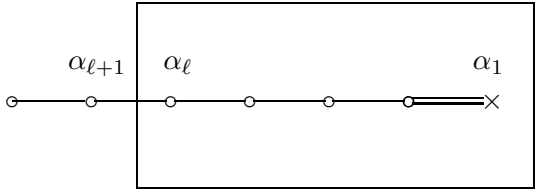
\vs
\begin{picture} (200, 50)

\put(24, 15){$\alpha_{\ell+1}$} \put(60,15){$\alpha_{\ell}$}

\put(0,0) {$\circ$}  \put(3,3){\line(1,0){27}} \put(30,0){$\circ$}
\put(34, 3){\line(1,0){27}}

\put(60,0) {$\circ$}  \put(64,3){\line(1,0){27}}

\put(90,0) {$\circ$}  \put(94,3){\line(1,0){27}}

\put(120,0) {$\circ$}  \put(124,3){\line(1,0){27}}

\put(150,0) {$\circ$}  \put(154,4){\line(1,0){27}}

\put(150,0) {$\circ$}  \put(154,2){\line(1,0){27}}

\put (180, 0) {$\times$}

\put(177, 15){$\alpha_{1}$}

\put(50, 40){\line(1,0){150}} \put(50,40){\line(0,-1){70}}

\put(50, -30){\line(1,0){150}} \put(200, -30){\line(0,1){70}}

\end{picture}

\vs \vs \vs \vs \vs

Then $\Delta(\tau)=\Delta(w)-\{\theta_{\delta}\} \cap \{
\theta_{\tau}\}$ and $\theta_{\tau}=\alpha_1 + \cdots +
\alpha_{\ell+1}$ and $\theta_{\delta}=\alpha_1 + \beta_1+
\beta_2$, where $\beta_1=\beta_2=\alpha_2+ \cdots +
\alpha_{\ell}$.
\end{pf}

\subsection{Conclusion} Recall that $k_{G/P}$ is defined by the
length of the chain with extremal node $\gamma$ and the branch
point(Section \ref{kGP}).

\begin{thm} \label{thm1} Let $G/P$ be a Hermitian symmetric space $\not= Q^{2n-1}$.
Let $\delta$ be a subgraph of $\CD(G)$ of type $A_{k}$ with
$\gamma$ an extremal  node and $X_w$ be the corresponding Schubert
variety. Suppose that  $k >k_{G/P}$. Then a subvariety $X$ of
$G/P$ with the homology class  $[X]=r[X_w]$ is contained in a
maximal linear space  in $G/P$. In particular, the cycle space
$\mathcal Z_k(G/P, [X_w])$ consists of Schubert varieties of type
$w$.

If $\delta$ is a maximal chain, then $X_w$ is Schur rigid.

\end{thm}

\begin{pf} By Proposition \ref{SchurHSS}, any integral variety of
$\CR_w$ is an integral variety of $\CB_w$ if $k > k_{G/P}$. By
Proposition \ref{SchubertHSS}, an integral variety of $\CB_w$ is
contained in a maximal linear space if $k > k_{G/P}$. Since a
subvariety of a linear space of degree 1 is again a linear space,
the cycle space $\mathcal Z_k(G/P, [X_w])$ consists of Schubert
varieties of type $w$.
\end{pf}

\begin{thm} \label{thm2} Let $G/P$ be a Hermitian symmetric space
associated to the simple root $\gamma$. Let $\delta$ be a subgraph
of $\CD(G)$, either

(1) of type $A_k$ with $\gamma$ non extremal node or

(2) of type $D_{\ell}$ with $\gamma$ an extremal node or

(3) of type $C_{\ell}$ with $\gamma=\alpha_{\ell}$.

\noindent  Then the corresponding Schubert variety $X_w$ is Schur
rigid.
\end{thm}

\begin{pf} By Proposition \ref{equivalence}, it suffices to prove
the Schubert rigidity and the equality $B_w=R_w$. The former is
given by Proposition \ref{smoothSchubert} and the latter is given
by Proposition \ref{smoothSchur} for the case (1) and by
Proposition \ref{SchurHSSD} for the case (2) and (3).
\end{pf}

By Proposition \ref{smoothdynkin}, Schubert varieties considered
in  Theorem \ref{thm1} and Theorem  \ref{thm2} are all the smooth
Schubert varieties in Hermitian symmetric spaces. So  we prove

\vs {\bf Main Theorem.} {\it Let $G/P$ be a Hermitian symmetric
space other than an odd dimensional quadric. Then any smooth
Schubert variety $X_w$ in $G/P$ is Schur rigid except when $X_w$
is a non-maximal linear space in $G/P$. Here, we consider $G/P$ as
a projective variety by the minimal equivariant embedding $G/P
\subset \PP(V)$. }

 \vs

\vs\vs

Research Institute of Mathematics

Seoul National University

San 56-1 Sinrim-dong Kwanak-gu

Seoul, 151-747 Korea

jhhong@math.snu.ac.kr


\begin{thebibliography}{99}


\bibitem[B]{B} R. Bryant, Rigidity and quasi-rigidity of extremal
cycles in compact Hermitian symmetric spaces, To appear in Annals
of mathematics studies 153, Princeton University Press.

\bibitem[BE]{BE} R. Baston and M. Eastwood, The Penrose Transform
Its interaction with representation theory, Oxford Science
Publication, 1989


\bibitem[BP]{BP}M. Brion, P. Polo, \emph{Generic singularities of
certain Schubert varieties}, Math. Z., 231, 301-324 (1999)

\bibitem[CH]{CH} I. Choe and J. Hong, \emph{Integral varieties of the
canonical cone structure on $G/P$}, Math. Ann. 329, 629-652 (2004)


\bibitem[F]{F} W. Fulton, Young tableaux with application to representatin theory and geometry,
London mathematical Society, Student Texts 35, Cambridge
University press, 1997

\bibitem[FH]{FH} W. Fulton and J. Harris, Representation Theory; A
First Course, Springer-Verlag, 1991


\bibitem[G]{G} A. B. Goncharov, \emph{Generalized conformal structures
on manifolds}, Selecta Mathematica Sovietica, vol 6, No 4, 1987.



\bibitem[K1]{K1} B. Kostant, \emph{Lie algebra cohomology and the
generalized Borel-weil theorem}, Ann. of Math. (2) 74 (1961), No.
2, 329-387.

\bibitem[K2]{K2} B. Kostant, \emph{Lie algebra cohomology and
generalized Schubert cells}, Ann. of Math. (2) 77 (1963), No. 1,
72-144.

\bibitem[LW]{LW} V. Lakshmibai and J. Weyman, \emph{Multiplicities
of Points on a Schubert variety in a Minuscule $G/P$}, Advances in
Mathematics 84, 179-208 (1990)



\bibitem[W]{W} M. Walters, \emph{Geometry and uniqueness of some
extreme subvarieties in complex Grassmannians}, Ph.D. thesis,
University of Michigan, 1997.


\end{thebibliography}
\end{document}